\documentclass[12pt]{article}
\usepackage{amsfonts}
\usepackage{amsmath}
\usepackage{hyperref}
\usepackage[utf8]{inputenc}
\usepackage{ stmaryrd }
\usepackage{multicol}
\usepackage{pgf,tikz}

\def\squarebox#1{\hbox to #1{\hfill\vbox to #1{\vfill}}}
\newcommand{\qed}{\hspace*{\fill}
\vbox{\hrule\hbox{\vrule\squarebox{.667em}\vrule}\hrule}\smallskip}
\newtheorem{teo}{Theorem}[section]
\newtheorem{lema}[teo]{Lemma}
\newtheorem{prop}[teo]{Proposition}
\newtheorem{cor}[teo]{Corollary}
\newtheorem{dfn}[teo]{Definition}
\newenvironment{demons}{\noindent {\bf Proof:}}{\hfill $\qed $ \newline}
\newenvironment{obs}{\noindent {\bf Remark:}}{\newline}
\newenvironment{exem}{\noindent {\bf Example:}}{\newline}
\newenvironment{nota}{\noindent {\bf Notation:}}{\newline}

\DeclareMathOperator{\Ad}{Ad}
\DeclareMathOperator{\ad}{ad}
\DeclareMathOperator{\Nij}{Nij}

\DeclareMathOperator{\tr}{tr}
\DeclareMathOperator{\hht}{ht}

\begin{document}

\title{Invariant Generalized Complex Structures on Partial Flag Manifolds}
\author{Carlos A. B. Varea\thanks{This work was supported by FAPESP grant 2016/07029-2. It was also financed in party by the Coordena\c{c}\~ao de Aperfei\c{c}oamento de Pessoal de Nível Superior - Brasil (CAPES) - Finance Code 001.}}
\date{}
\maketitle

\begin{abstract}
The aim of this paper is to classify all invariant generalized complex structure on a partial flag manifold $\mathbb{F}_\Theta$ with at most four isotropy summands. To classify them all we proved that an invariant generalized almost complex structure on $\mathbb{F}_\Theta$ is `constant' in each component of the isotropy representation.
\end{abstract}
\noindent
\\%

\noindent%

\section{Introduction}

The subject matter of this paper are invariant generalized complex structures on partial flag manifolds of semisimple Lie groups. A generalized complex structure is a differential geometric structure introduced by Hitchin in \cite{H} and further developed by Gualtieri \cite{G1}, with the purpose of studying complex and symplectic structures in a unique framework. We refer to Guatieri \cite{G1}, \cite{G2} and Cavalcanti \cite{C} for the foundations of the theory of generalized complex structures.

In this paper we consider partial flag manifolds of complex Lie groups, sometimes it is also called generalized flag manifolds. Let $\mathfrak{g}$ be a complex semisimple Lie algebra, $G$ be a connected Lie group with Lie algebra $\mathfrak{g}$. Then a partial flag manifold is the homogeneous space $\mathbb{F}_\Theta = G/P_\Theta$ where $P_\Theta$ is a parabolic subgroup associated to the parabolic subalgebra $\mathfrak{p}_\Theta$. If $U$ is a compact real form of $G$ then $U$ acts transitively on $\mathbb{F}_\Theta$ so that we have also the homogeneous space $\mathbb{F}_\Theta = U/K$ where $K = P_\Theta \cap U$ is the centralizer of a torus. We are concerned with $U$-invariant generalized complex structures on $\mathbb{F}_\Theta$.

Our approach to study invariant generalized complex structures on $\mathbb{F}_\Theta$ is to reduce the problem at the origin. In \cite{VS} the authors give a description of all invariant generalized complex structures on a maximal flag manifold, that is, when $\Theta = \emptyset$. It is done restricting the invariant generalized complex structure $\mathcal{J}$ to the subspace $\mathfrak{u}_\alpha \oplus \mathfrak{u}_\alpha ^\ast$ for each positive root $\alpha$, we denote such restriction by $\mathcal{J}_\alpha$. They proved that there are only two types of generalized complex structure on $\mathfrak{u}_\alpha \oplus \mathfrak{u}_\alpha ^\ast$, complex and noncomplex type. 

There exists a complete classification of all partial flag manifolds with two, three and four isotropy summands, see \cite{AC1}, \cite{Kim} and \cite{AC2}. Based on this classification our problem reduce to study what happen in an irreducible component of the isotropy representation. It is known that an invariant complex structure on a partial flag manifold can be described by signs, that is, if $J$ is an invariant complex structure on a flag manifold $\mathbb{F}_\Theta$, then $J(X_\alpha) = \varepsilon_\alpha i X_\alpha$ where $\varepsilon_\alpha = \pm 1$. Moreover, if $\mathfrak{m}_j$ is an irreducible component of the isotropy representation, it is known that $\varepsilon_\alpha$ is constant inside $\mathfrak{m}_j$. Analogously, it is also known, that an invariant metric on $\mathbb{F}_\Theta$ is constant inside an irreducible component of the isotropy representation. Motivated by these results it is reasonable to expect that, somehow, an invariant generalized complex structure $\mathcal{J}$ on a partial flag manifold $\mathbb{F}_\Theta$ is `constant' when restricted to $\mathfrak{m}_j \oplus \mathfrak{m}_j ^\ast$, where $\mathfrak{m}_j$ is an irreducible component of the isotropy representation.

In this approach the first step is to determine the behaviour of an invariant generalized complex structure $\mathcal{J}$ on a partial flag manifold $\mathbb{F}_\Theta$ restricted to $\mathfrak{m}_j \oplus \mathfrak{m}_j ^\ast$, where $\mathfrak{m}_j$ is an irreducible component of the isotropy representation. Relying on the invariance of $\mathcal{J}$ we proved that the restriction $\mathcal{J}|_{\mathfrak{m}_j \oplus \mathfrak{m}_j ^\ast}$ has constant type, that is, if $\alpha$ is a positive root in the sum $\mathfrak{m}_j = \sum_\beta  \mathfrak{u}_\beta$ such that $\mathcal{J}_\alpha$ is of (non)complex type, then $\mathcal{J}_\beta$ is also of non(complex) type for all $\beta$ appearing in the sum $\mathfrak{m}_j = \sum_\beta \mathfrak{u}_\beta$. The proof of this fact has an independent interest of the classification presented in this paper.

Once we have the description of the invariant generalized complex structures on a partial flag manifold and we have the classification of all partial flag manifolds with at most four isotropy summands, we proceed to analyze the integrability of such structures. In \cite{VS} the authors proved that the integrability of an invariant generalized complex structure depends on triples of roots $(\alpha,\beta,\alpha+\beta)$ by means of analyzing the Nijenhuis operator restricted to the $i$-eigenbundle of the triple $(\mathcal{J}_\alpha, \mathcal{J}_\beta, \mathcal{J}_{\alpha+\beta})$.

The paper is organized as follows: In Section \ref{GFM} we study the structure of a partial flag manifold $\mathbb{F}_\Theta = G/P_\Theta$ of a semisimple Lie group $G$. In particular, we present a way to describe the components of the isotropy representation. In Section \ref{GCG} we recall some basic concepts of generalized complex geometry. We also present the classification of the invariant generalized complex structures on maximal flag manifolds in order to describe such structures on partial flag manifolds. In the end of this section we prove that an invariant generalized complex structure on $\mathbb{F}_\Theta$ is `constant' inside the components of the isotropy representation. In Section \ref{CAP4} we present the classification of all flag manifolds with two, three and four isotropy summands and we classify all invariant integrable generalized complex structures on these flag manifolds.

\section{Partial Flag Manifolds}\label{GFM}

Let $\mathfrak{g}$ be a semisimple Lie algebra and $G$ be a connected Lie group with Lie algebra $\mathfrak{g}$. Given $\mathfrak{h}$ a Cartan subalgebra of $\mathfrak{g}$, let $\Pi$ be a root system of $\mathfrak{g}$ relative to $\mathfrak{h}$. We can decompose $\mathfrak{g}$ as
\[
\mathfrak{g} = \mathfrak{h} \oplus \sum_{\alpha \in \Pi} \mathfrak{g}_\alpha,
\]
where $\mathfrak{g}_\alpha = \{ X\in \mathfrak{g} \ | \ [H,X] = \alpha (H)X \}$ is the root space associated to $\alpha$. The Cartan--Killing form $\langle X,Y \rangle = \tr (\ad (X) \ad (Y) )$ of $\mathfrak{g}$ is nondegenerate on $\mathfrak{h}$. Given $\alpha \in \mathfrak{h}^\ast$ we let $H_\alpha$ be defined by $\alpha ( \cdot ) = \langle H_\alpha , \cdot \rangle$, and denote by $\mathfrak{h}_\mathbb{R}$ the real subspace generated by $H_\alpha$, $\alpha \in \Pi$.

Let $\Pi^+ \subset \Pi$ be a choice of positive roots and denote by $\Sigma$ the corresponding simple root system. 

\begin{dfn}
The Borel subalgebra $\mathfrak{b}$ of $\mathfrak{g}$ is the maximal solvable subalgebra defined by
\[
\mathfrak{b} = \mathfrak{h} \oplus \sum_{\alpha \in \Pi^+} \mathfrak{g}_\alpha.
\]
We will say that a subalgebra $\mathfrak{p}$ of $\mathfrak{g}$ is a parabolic subalgebra when it contains the Borel subalgebra, that is, $\mathfrak{b} \subset \mathfrak{p}$.
\end{dfn}

Let $\Theta$ be a subset of $\Sigma$ and denote by $\langle \Theta \rangle$ the set of roots generated by $\Theta$, that is, if $\Theta = \{ \alpha_{i_1}, \cdots , \alpha_{i_k} \}$ and $\alpha$ is a root such that $\alpha \in \langle \Theta \rangle$ then $\alpha = \sum_{i=1} ^k n_i \alpha_{j_1}$ for some $n_i \in \mathbb{Z}$. Let $\langle \Theta \rangle ^+ = \Pi ^+ \cap \langle \Theta \rangle$, then we have that
\[
\mathfrak{p}_\Theta = \mathfrak{h} \oplus \sum_{\alpha \in \langle \Theta \rangle ^+} \mathfrak{g}_{-\alpha} \oplus \sum_{\alpha \in \Pi^+} \mathfrak{g}_\alpha
\]
is a parabolic subalgebra of $\mathfrak{g}$, since it contains the Borel subalgebra $\mathfrak{b} \subset \mathfrak{p}_\Theta$. 

\begin{dfn}
The partial flag manifold $\mathbb{F}_\Theta$ associated to $\mathfrak{p}_\Theta$ is the homogeneous space $\mathbb{F}_\Theta = G/P_\Theta$, where $P_\Theta$ is the parabolic subgroup generated by $\mathfrak{p}_\Theta$.
\end{dfn}

In particular, if $\Theta = \emptyset$ we have $\mathfrak{p}_\Theta = \mathfrak{b}$. In this case $\mathbb{F}_\Theta = \mathbb{F}$ is called maximal flag manifold.

Let $\mathfrak{u}$ be a compact real form of $\mathfrak{g}$, to know, the real subalgebra 
\[
\mathfrak{u} = \textnormal{span}_\mathbb{R} \{ i\mathfrak{h}_\mathbb{R},A_\alpha, S_\alpha \ : \ \alpha\in \Pi^+ \}
\]
where $A_\alpha = X_\alpha - X_{-\alpha}$ and $S_\alpha = i(X_\alpha + X_{-\alpha})$. Denote by $U = \exp \mathfrak{u}$ the correspondent compact real form of $G$. Then the real representation $\mathbb{F}_\Theta = U/K$ is obtained by the transitive action of $U$ on $G/P_\Theta$, where the closed connected subgroup $K = P_\Theta \cap U$ can be identified with the centralizer $C(T)$ of a torus $T\subset U$. We have that $\mathfrak{k} = \mathfrak{p}_\Theta \cap \mathfrak{u}$ is the Lie algebra of $K$ and set $\mathfrak{m}$ the orthogonal complement of $\mathfrak{k}$ on $\mathfrak{u}$, that is, $\mathfrak{u} = \mathfrak{k} \oplus \mathfrak{m}$. Thus, the tangent space of $\mathbb{F}_\Theta$ at the origin can be identified with $\mathfrak{m}$. 

Since $\mathbb{F}_\Theta$ is a reductive homogeneous space, the isotropy representation of $K$ on the tangent space of $\mathbb{F}_\Theta$ at the origin is equivalent to the adjoint representation $\Ad|_K$ restricted to $\mathfrak{m}$. It is known that this representation is completely reductive, so we can decompose $\mathfrak{m}$ into mutually non equivalent irreducible $\Ad (K)$-module as 
\[
\mathfrak{m} = \mathfrak{m}_1 \oplus \cdots \oplus \mathfrak{m}_s,
\]
where $\mathfrak{m}_i$ is an irreducible component of the isotropy representation, for each $i= 1,\cdots, s$.

Let $\Sigma = \{ \alpha_1, \cdots, \alpha_l\}$ and $\Theta = \{\alpha_{i_1},\cdots , \alpha_{i_k} \} \subset \Sigma$, then we can write $\Sigma \backslash \Theta = \{ \alpha_{j_1}, \cdots ,\alpha_{j_r}\}$ where $k + r = l$. For integers $s_1, \cdots , s_r$ with $(s_1,\cdots ,s_r) \not= (0,\cdots ,0)$ we set
\[
\Pi ^\mathfrak{m} (s_1,\cdots ,s_r) = \left\lbrace \sum_{i=1} ^l n_i\alpha_i \in \Pi^+ \ : \ n_{j_1} = s_1, \cdots , n_{j_r} = s_r \right\rbrace . 
\]
Then $\Pi_\mathfrak{m}^+ = \Pi^+ \backslash \langle \Theta \rangle = \displaystyle\bigcup_{s_1,\cdots ,s_r} \Pi^\mathfrak{m} (s_1,\cdots ,s_r)$. It was proved in \cite{ACS} that for $\Pi^\mathfrak{m} (s_1,\cdots ,s_r) \not= \emptyset$ we define an $\Ad(K)$-invariant subspace $\mathfrak{m}(s_1,\cdots,s_r)$ of $\mathfrak{u}$ by
\[
\mathfrak{m}(s_1,\cdots,s_r) = \sum_{\alpha \in \Pi^\mathfrak{m} (s_1,\cdots ,s_r)} \mathfrak{u}_\alpha
\]
where $\mathfrak{u}_\alpha = \textnormal{span}_\mathbb{R} \{ A_\alpha, S_\alpha \}$. Moreover, they proved  that we have a decomposition of $\mathfrak{m}$ into mutually non equivalent irreducible $\Ad(K)$-module given by
\[
\mathfrak{m} = \sum _{s_1,\cdots,s_r} \mathfrak{m}(s_1,\cdots ,s_r)
\]

\begin{exem}
Let $G$ be a Lie group with Lie algebra $\mathfrak{g}$ of type $B_3$. Then we have $\Sigma = \{\alpha_1,\alpha_2,\alpha_3\}$ a simple root system and $\Pi^+ = \{ \alpha_1,\alpha_2,\alpha_3,\alpha_1+\alpha_2,\alpha_2+\alpha_3, \alpha_1+\alpha_2+\alpha_3,\alpha_2+2\alpha_3,\alpha_1+\alpha_2+2\alpha_3,\alpha_1+2\alpha_2+2\alpha_3\}$ the set of positive roots. Let $\Theta = \{ \alpha_3\}$ and let $\mathbb{F}_\Theta$ be the corresponding flag manifold. 

We have that $\mathfrak{m}$ decomposes into four non equivalent irreducible components, that is, $\mathfrak{m} = \mathfrak{m}_1 \oplus \mathfrak{m}_2 \oplus \mathfrak{m}_3 \oplus \mathfrak{m}_4$. In fact, following the notation fixed above, we have $\Pi ^\mathfrak{m}(1,0) = \{\alpha_1 \}$, $\Pi ^\mathfrak{m} (0,1) = \{ \alpha_2,\alpha_2+\alpha_3, \alpha_2+2\alpha_3\}$, $\Pi^\mathfrak{m} (1,1) = \{ \alpha_1+\alpha_2, \alpha_1+\alpha_2+\alpha_3, \alpha_1+\alpha_2+2\alpha_3\}$ and $\Pi^\mathfrak{m} (1,2) = \{ \alpha_1 +2\alpha_2+2\alpha_3\}$. Therefore $\mathfrak{m}_1 = \mathfrak{m}(1,0) = \mathfrak{u}_{\alpha_1}$, $\mathfrak{m}_2 = \mathfrak{m}(0,1) = \mathfrak{u}_{\alpha_2} \oplus \mathfrak{u}_{\alpha_2+\alpha_3} \oplus \mathfrak{u}_{\alpha_2+2\alpha_3}$, $\mathfrak{m}_3 = \mathfrak{m}(1,1) = \mathfrak{u}_{\alpha_1+\alpha_2} \oplus \mathfrak{u}_{\alpha_1+\alpha_2+\alpha_3} \oplus \mathfrak{u}_{\alpha_1+\alpha_2+2\alpha_3}$ and $\mathfrak{m}_4 = \mathfrak{m}(1,2) = \mathfrak{u}_{\alpha_1+2\alpha_2+2\alpha_3}$.
\end{exem}

Now, for each $\alpha_i \in \Sigma$ consider its dual $H_i$, that is, $\alpha_i (H_j) = \delta_{ij}$. Let $H_0 = a_1 H_{j_1} + \cdots + a_r H_{j_r}$ be a generic element of the positive Weyl chamber. It is clear that $\alpha (H_0) = 0$ for all $\alpha \in \Pi \backslash \langle \Theta \rangle$ and $[H_0, X] = \alpha (H_0)X$ for all $X\in \mathfrak{m}$. Then we can describe the irreducible components of $\mathfrak{m}$ by the nonzero eigenvalues of $\ad (H_0)$. That is, given $\alpha \in \Pi^+$ where $\alpha = \sum_{i=1} ^l n_i\alpha_i$, then we have that $\ad(H_0)X = \alpha (H_0) X = (n_{j_1} a_1 + \cdots + n_{j_r}a_r)X$ for $X\in \mathfrak{u}_\alpha$. 


Observe that
\begin{eqnarray*}
\Pi ^\mathfrak{m} (s_1,\cdots ,s_r) & = & \left\lbrace \alpha = \sum_{i=1} ^l n_i\alpha_i \in \Pi^+ \ : \ n_{j_1} = s_1, \cdots , n_{j_r} = s_r \right\rbrace \\
& = & \left\lbrace \alpha \in \Pi^+ \ : \alpha (H_0) = a_1s_1 + \cdots a_rs_r \right\rbrace .
\end{eqnarray*}
Therefore
\[
\mathfrak{m} = \sum _{s_1,\cdots,s_r} \mathfrak{m}(s_1,\cdots ,s_r)
\]
is a decomposition of $\mathfrak{m}$ into mutually non equivalent irreducible $\Ad(K)$-module, where each $\mathfrak{m}(s_1,\cdots,s_r)$ is related to an eigenvalue of $\ad (H_0)$.\\

\begin{exem}
Let $G$ be a Lie group with Lie algebra $\mathfrak{g}$ of type $B_3$. Then we have $\Sigma = \{\alpha_1,\alpha_2,\alpha_3\}$ a simple root system and $\Pi^+ = \{ \alpha_1,\alpha_2,\alpha_3,\alpha_1+\alpha_2,\alpha_2+\alpha_3, \alpha_1+\alpha_2+\alpha_3,\alpha_2+2\alpha_3,\alpha_1+\alpha_2+2\alpha_3,\alpha_1+2\alpha_2+2\alpha_3\}$ the set of positive roots. Let $\Theta = \{ \alpha_3\}$ and let $\mathbb{F}_\Theta$ be the corresponding flag manifold. Let $H_0 = a_1\alpha_1 + a_2\alpha_2$ be a generic element of the positive Weyl chamber. 

It is easy to prove that $\alpha_1 (H_0) = a_1$, $\alpha_2 (H_0) = (\alpha_2+\alpha_3)(H_0) = (\alpha_2+2\alpha_3)(H_0) = a_2$, $(\alpha_1+\alpha_2)(H_0) = ( \alpha_1+\alpha_2+\alpha_3)(H_0) = (\alpha_1+\alpha_2+2\alpha_3)(H_0) = a_1+a_2$ and $(\alpha_1 + 2\alpha_2 +2\alpha_3)(H_0) = a_1+2a_2$. Then the irreducible components are given by $\mathfrak{m}_1 = \mathfrak{m}(1,0) = \mathfrak{u}_{\alpha_1}$, $\mathfrak{m}_2 = \mathfrak{m}(0,1) = \mathfrak{u}_{\alpha_2} \oplus \mathfrak{u}_{\alpha_2+\alpha_3} \oplus \mathfrak{u}_{\alpha_2+2\alpha_3}$, $\mathfrak{m}_3 = \mathfrak{m}(1,1) = \mathfrak{u}_{\alpha_1+\alpha_2} \oplus \mathfrak{u}_{\alpha_1+\alpha_2+\alpha_3} \oplus \mathfrak{u}_{\alpha_1+\alpha_2+2\alpha_3}$ and $\mathfrak{m}_4 = \mathfrak{m}(1,2) = \mathfrak{u}_{\alpha_1+2\alpha_2+2\alpha_3}$.
\end{exem}

Summing up, to get the decomposition of $\mathfrak{m}$ into mutually non equivalent irreducible $\Ad (K)$-module is enough to look for the eigenvalues of $\ad(H_0)$, where $H_0$ is a generic element of the positive Weyl chamber. 

\section{Generalized Complex Geometry}\label{GCG}
In this section we introduce the basic definitions of generalized complex
geometry. For more details see Gualtieri \cite{G1}. Let $M$ be a
smooth $n$-dimensional manifold, then the sum of the tangent and cotangent
bundle $TM\oplus T^{\ast }M$ is endowed with a natural symmetric bilinear
form with signature $(n,n)$ defined by 
\begin{equation*}
\langle X+\xi ,Y+\eta \rangle =\frac{1}{2}(\xi (Y)+\eta (X)).
\end{equation*}%
Furthermore, the Courant bracket is a skew-symmetric bracket defined on
smooth sections of $TM\oplus T^{\ast }M$ by 
\begin{equation*}
\lbrack X+\xi ,Y+\eta ]=[X,Y]+\mathcal{L}_{X}\eta -\mathcal{L}_{Y}\xi -\frac{%
1}{2}d\left( i_{X}\eta -i_{Y}\xi \right) .
\end{equation*}

\begin{dfn}
A generalized almost complex structure on $M$ is a map $\mathcal{J}\colon
TM\oplus T^{\ast }M\rightarrow TM\oplus T^{\ast }M$ such that $\mathcal{J}%
^{2}=-1$ and $\mathcal{J}$ is an isometry of the bilinear form $\langle
\cdot ,\cdot \rangle $. The generalized almost complex structure $\mathcal{J}
$ is said to be integrable to a generalized complex structure when its $i$%
-eigenbundle $L\subset (TM\oplus T^{\ast }M)\otimes \mathbb{C}$ is Courant
involutive.
\end{dfn}

Note that, given $L$ a maximal isotropic sub-bundle of $TM\oplus T^{\ast }M$
(or its complexification) then $L$ is Courant involutive if and only if $\Nij%
|_{L}=0$, where $\Nij$ is the Nijenhuis operator defined by 
\begin{equation}\label{GeneralNij}
\Nij(A,B,C)=\frac{1}{3}\left( \langle \lbrack A,B],C\rangle +\langle \lbrack
B,C],A\rangle +\langle \lbrack C,A],B\rangle \right) .
\end{equation}

\begin{exem}
The basic examples of generalized complex structures come from complex
and symplectic structures. If $J$ and $\omega $ are complex and symplectic
structures respectively on $M$, then 
\begin{equation*}
\mathcal{J}_{J}=\left( 
\begin{array}{cc}
-J & 0 \\ 
0 & J^{\ast }%
\end{array}%
\right) \ \mathrm{and}\ \mathcal{J}_{\omega }=\left( 
\begin{array}{cc}
0 & -\omega ^{-1} \\ 
\omega & 0%
\end{array}%
\right)
\end{equation*}%
are generalized complex structures on $M$.
\end{exem}

In \cite{VS} is presented a description of the invariant generalized complex structures on maximal flag manifolds, aiming to describe this structures we identify $\mathfrak{u}\cong \mathfrak{u}^\ast$ and $\mathfrak{u}_\alpha\cong \mathfrak{u}_\alpha^\ast$ by means of the Kirillov--Kostant--Souriau (KKS) symplectic form on $\textnormal{Ad}(U)(H)$ which at the origin $b_0$ is given by 
$$\omega_{b_0}(\widetilde{X},\widetilde{Y})=\langle H,[X,Y]\rangle,$$
for all $X,Y\in\mathfrak{u}$ and where $\widetilde{X}=\textnormal{ad}(X)$ denotes de fundamental vector field associated to the adjoint action. The elements of $\mathfrak{u}_\alpha^\ast$ will be denoted by $A_\alpha^\ast$ and $S_\alpha^\ast$ where
$$X^\ast=\dfrac{1}{\langle H,H_\alpha\rangle}\omega(\widetilde{X},\cdot).$$
For every $\alpha\in \Pi^+$, we will write all our object restrict to $\mathfrak{u}_\alpha\oplus\mathfrak{u}_\alpha^\ast$ in terms of the basis $\lbrace A_\alpha,S_\alpha,-S_\alpha^\ast,A_\alpha^\ast\rbrace$. As a first consequence of these identifications we get that for triples of roots $\alpha,\beta,\gamma\in\Pi^+$, the Nijenhuis tensor \eqref{GeneralNij} can be rewritten as
\begin{equation}\label{ParticularNij}
N(A,B,C)=\dfrac{1}{2}(k_\gamma\langle H,[C_2,[A_1,B_1]] \rangle+k_\alpha\langle H,[A_2,[B_1,C_1]] \rangle+k_\beta\langle H,[B_2,[C_1,A_1]] \rangle),
\end{equation}
where $A=A_1+A_2^\ast\in \mathfrak{u}_\alpha\oplus \mathfrak{u}_\alpha^\ast$, $B=B_1+B_2^\ast\in\mathfrak{u}_\beta\oplus\mathfrak{u}_\beta^\ast$, $C=C_1+C_2^\ast\in\mathfrak{u}_\gamma\oplus \mathfrak{u}_\gamma^\ast$, and $k_\alpha=\dfrac{1}{\langle H,H_\alpha\rangle}$ for all $\alpha\in \Pi^+$.

The invariant generalized almost complex structures $\mathcal{J}$ on $\mathbb{F}$ are those whose restriction to $\mathfrak{u}_\alpha\oplus\mathfrak{u}_\alpha^\ast$ have the form either
$$ \mathcal{J}_\alpha = \pm \mathcal{J}_0=\pm \left( 
\begin{array}{cccc}
0 & -1 & 0 & 0\\ 
1 & 0 & 0 & 0\\
0 & 0 & 0 & -1\\
0 & 0 & 1 & 0
\end{array}%
\right):\qquad \textbf{complex type},$$
or else
$$\mathcal{J}_\alpha=\left( 
\begin{array}{cccc}
a_\alpha & 0 & 0 & -x_\alpha\\ 
0 &  a_\alpha& x_\alpha & 0\\
0 & -y_\alpha & -a_\alpha & 0\\
y_\alpha & 0 & 0 & -a_\alpha
\end{array}%
\right):\qquad \textbf{noncomplex type},$$
with $a_\alpha,x_\alpha,y_\alpha\in\mathbb{R}$ such that $a_\alpha^2=x_\alpha y_\alpha-1$.\\

\begin{nota}
Let $\mathcal{J}$ be an invariant almost generalized complex structure on a flag manifold $\mathbb{F}$. We will denote by $\mathcal{J}_\alpha$ the restriction of $\mathcal{J}$ to $\mathfrak{u}_\alpha \oplus \mathfrak{u}_\alpha ^\ast$. 
\end{nota}

The integrability of the invariant generalized almost complex structures given above depends on analyzing what happens with triples of the form $(\mathcal{J}_\alpha,\mathcal{J}_\beta,\mathcal{J}_{\alpha+\beta})$
associated to triples of positive roots $(\alpha,\beta,\alpha+\beta)$. 
Accordingly, we obtain that $\mathcal{J}$
 is integrable if and only if for each triple of positive roots we have that $(\mathcal{J}_\alpha,\mathcal{J}_\beta,\mathcal{J}_{\alpha+\beta})$ 
corresponds to one of the rows of the following table 

\begin{table}[htb!]
\begin{center}
\begin{tabular}{|c|c|c|}
\hline $\mathcal{J}_\alpha$ & $\mathcal{J}_\beta$ & $\mathcal{J}_{\alpha + \beta}$ \\
\hline
complex ($\mathcal{J}_\alpha = \pm \mathcal{J}_0$) & complex ($\mathcal{J}_\beta = \pm \mathcal{J}_0$) & complex ($\mathcal{J}_{\alpha+\beta} = \pm \mathcal{J}_0$)\\
complex ($\mathcal{J}_\alpha = \pm \mathcal{J}_0$) & complex ($ \mathcal{J}_\beta = \mp \mathcal{J}_0$) & complex ($\mathcal{J}_{\alpha+\beta} = \pm \mathcal{J}_0$)\\
complex ($\mathcal{J}_\alpha = \pm \mathcal{J}_0$) & complex ($ \mathcal{J}_\beta = \mp \mathcal{J}_0$) & complex ($\mathcal{J}_{\alpha+\beta} = \mp \mathcal{J}_0$)\\
noncomplex  & complex ($ \mathcal{J}_\beta = \pm \mathcal{J}_0 $) & complex ($\mathcal{J}_{\alpha+\beta}= \pm \mathcal{J}_0)$ \\
complex ($\mathcal{J}_\alpha =  \pm \mathcal{J}_0$) & noncomplex & complex $(\mathcal{J}_{\alpha+\beta} = \pm \mathcal{J}_0$) \\
complex $(\mathcal{J}_\alpha = \pm \mathcal{J}_0)$ & complex ($\mathcal{J}_\beta = \mp \mathcal{J}_0)$ & noncomplex\\
noncomplex & noncomplex & noncomplex \\ \hline
\end{tabular} 
\caption{Integrability conditions.}
\label{Jintegravel}
\end{center}
\end{table}
\noindent where, in last row, the following extra condition is required for integrability:
\begin{equation}\label{Integrabilityconditions}
\left\lbrace \begin{array}{cc}
a_{\alpha+\beta} x_\alpha x_\beta - a_\beta x_\alpha x_{\alpha+\beta} - a_\alpha x_\beta x_{\alpha+\beta} = 0 \\
x_\alpha x_\beta - x_\alpha x_{\alpha+\beta} - x_\beta x_{\alpha+\beta} = 0.
\end{array}\right.
\end{equation} 
 
These structures are invariant by the torus action, then when we consider an invariant generalized almost complex structure $\mathcal{J}$ on a partial flag manifold it must be invariant by the action of a centralizer of a torus then, in particular, it is invariant by the torus action. Therefore, the statement done above for a maximal flag manifold still hold for partial flag manifolds. We must to be careful with the integrability conditions, because in a maximal flag manifold the irreducible components of $\mathfrak{m}$ are 1-dimensional, what is not true when we consider partial flag manifolds. In some sense, it is expected $\mathcal{J}$ to be `constant' on the irreducible components, that is, let $\mathfrak{m} = \mathfrak{m}_1 \oplus \cdots \oplus \mathfrak{m}_s$ then we expect that $\mathcal{J}$ restricted to $\mathfrak{m}_i \oplus \mathfrak{m}_i ^\ast$ be `constant'.

\begin{prop}
Let $\mathbb{F}_\Theta$ be a partial flag manifold and $\mathcal{J}$ be an invariant generalized almost complex structure on $\mathbb{F}_\Theta$. Suppose that $\mathfrak{m} = \mathfrak{m}_1 \oplus \cdots \oplus \mathfrak{m}_s$. Then $\mathcal{J}$ restricted to $\mathfrak{m}_i \oplus \mathfrak{m}_i ^\ast$ has constant type, that is, if $\alpha$ and $\beta$ are associated to $\mathfrak{m}_i$, then $\mathcal{J}_\alpha$ and $\mathcal{J}_\beta$ has the same type (complex or noncomplex type), for all $i=1,\cdots,s$. Moreover, $\mathcal{J}_\alpha = \mathcal{J}_\beta$. 
\end{prop}
\begin{demons}
Consider $\mathfrak{m}_k$ an irreducible component of $\mathfrak{m}$ and, abusing the notation, let $\alpha,\beta \in \mathfrak{m}_k$. It is well known that there is a root $\gamma \in \langle \Theta \rangle$ such that $\alpha + \gamma = \beta$. Consider $A_\gamma \in \mathfrak{u}_\gamma$, then
\begin{eqnarray*}
\ad (A_\gamma)A_\alpha = m_{\gamma,\alpha} A_\beta + m_{-\gamma,\alpha} A_{\gamma-\alpha} \\
\ad (A_\gamma)S_\alpha = m_{\gamma,\alpha} S_\beta + m_{\gamma,-\alpha} A_{\gamma-\alpha}.
\end{eqnarray*}
Thus, on the matrix of $(\ad \oplus \ad ^\ast)(A_\gamma)$ appears a block of the form
\[
M = \left( \begin{array}{cccc}
m_{\gamma,\alpha} & 0 & 0 & 0\\
0 & m_{\gamma,\alpha} & 0 & 0\\
0 & 0 & \ast & \ast \\
0 & 0 & \ast & \ast 
\end{array} \right).
\]
Since $\mathcal{J}$ is an invariant generalized almost complex structure, then $\mathcal{J}$ must commute with $(\ad \oplus \ad ^\ast)(A_\gamma)|_{\mathfrak{m}\oplus \mathfrak{m}^\ast}$, that is, we must have
\begin{equation}\label{jcomuta}
\mathcal{J}\circ (\ad \oplus \ad ^\ast)(A_\gamma)|_{\mathfrak{m}\oplus \mathfrak{m}^\ast} = (\ad \oplus \ad ^\ast)(A_\gamma)|_{\mathfrak{m}\oplus \mathfrak{m}^\ast} \circ \mathcal{J}.
\end{equation}
But, if equation (\ref{jcomuta}) is satisfied then we have
\begin{equation}\label{Mcomuta}
\mathcal{J}_\beta \cdot M = M \cdot \mathcal{J}_\alpha.
\end{equation}
From (\ref{Mcomuta}) we get that $\mathcal{J}_\alpha$ and $\mathcal{J}_\beta$ have the same type. Moreover, (\ref{Mcomuta}) ensure that $\mathcal{J}_\alpha = \mathcal{J}_\beta$.
\end{demons}

\begin{nota}
For simplicity, we will abuse the notation and write $\alpha \in \mathfrak{m}_i$ when $\alpha$ is a root that appears in the sum $\mathfrak{m}_i = \sum_\beta \mathfrak{u}_\beta$.
\end{nota}

Therefore, we conclude that an invariant generalized almost complex structure $\mathcal{J}$ on $\mathbb{F}_\Theta$ is `constant' when restricted to $\mathfrak{m}_k$. In this case, constant means that $\mathcal{J}$ have the same type on $\mathfrak{m}_k \oplus \mathfrak{m}_k ^\ast$, that is, $\mathcal{J}_\alpha$ has the same type for all $\alpha \in \mathfrak{m}_k$. More than that, we have $\mathcal{J}_\alpha = \mathcal{J}_\beta$ for all $\alpha,\beta \in \mathfrak{m}_k$.

\section{Invariant Generalized Complex Structures on $\mathbb{F}_\Theta$}\label{CAP4}

The aim of this section is to analyse the integrability conditions of an invariant generalized almost complex structure $\mathcal{J}$ on a partial flag manifold $\mathbb{F}_\Theta$. We know that the integrability condition is given analysing a triple of roots $(\alpha,\beta,\alpha+\beta)$. For each triple of roots $(\alpha,\beta,\alpha+\beta)$ we have two possibilities:
\begin{enumerate}
\item $\alpha,\beta \in \mathfrak{m}_i$ and $\alpha+\beta \in \mathfrak{m}_j$;

\item $\alpha \in \mathfrak{m}_i$, $\beta \in \mathfrak{m}_j$ and $\alpha+\beta \in \mathfrak{m}_k$.
\end{enumerate}

The first case is when $\alpha, \beta \in \mathfrak{m}_i$ and $\alpha+\beta \in \mathfrak{m}_j$, with $i\not= j$. If $\mathcal{J}$ is integrable, by Table \ref{Jintegravel}, we must have $\mathcal{J}|_{\mathfrak{m}_i \oplus \mathfrak{m}_i ^\ast}$ and $\mathcal{J}|_{\mathfrak{m}_j \oplus \mathfrak{m}_j ^\ast}$ of the same type. Moreover, if $\mathcal{J}|_{\mathfrak{m}_i \oplus \mathfrak{m}_i ^\ast}$ is of complex type with $\mathcal{J}_\alpha= \mathcal{J}_\beta = \mathcal{J}_0$ ($=-\mathcal{J}_0$) then $\mathcal{J}|_{\mathfrak{m}_j \oplus \mathfrak{m}_j ^\ast}$ must be of complex type with $\mathcal{J}_{\alpha+\beta} = \mathcal{J}_0$ ($=-\mathcal{J}_0$). 

\begin{prop}
Let $\mathcal{J}$ be an invariant generalized almost complex structure on $\mathbb{F}_\Theta$, where $\mathfrak{m} = \mathfrak{m}_1 \oplus \cdots \oplus \mathfrak{m}_s$. If $\mathcal{J}$ is integrable and there exists a triple of roots $(\alpha,\beta,\alpha+\beta)$ such that $\alpha,\beta \in \mathfrak{m}_i$ and $\alpha+\beta \in \mathfrak{m}_j$ with $i\not= j$, then $\mathcal{J}|_{\mathfrak{m}_i \oplus \mathfrak{m}_i ^\ast}$ and $\mathcal{J}|_{\mathfrak{m}_j \oplus \mathfrak{m}_j ^\ast}$ have the same type.
\end{prop}

In particular:

\begin{cor}\label{2somandos}
Let $\mathcal{J}$ be an invariant almost complex structure on $\mathbb{F}_\Theta$, where $\mathfrak{m} = \mathfrak{m}_1 \oplus \mathfrak{m}_2$. If $\mathcal{J}$ is integrable there exists a triple of roots $(\alpha,\beta,\alpha+\beta)$ such that $\alpha,\beta$ are in the same component and $\alpha+\beta$ is on the other component, then $J_\gamma$ is of the same type, complex or noncomplex type, for all $\gamma \in \mathfrak{m}$.
\end{cor}

\begin{obs}
Observe that if $\mathcal{J}$ is an invariant generalized complex structure on a flag manifold $\mathbb{F}_\Theta$ such that $\mathfrak{m}= \mathfrak{m}_1 \oplus \mathfrak{m}_2$, by Corollary \ref{2somandos}, we have that $\mathcal{J}_\alpha$ is of the same type for all root $\alpha$. If $\mathcal{J}_\alpha$ is of complex type, then we have that $\mathcal{J}_\alpha = \mathcal{J}_0$ or $\mathcal{J}_\alpha = - \mathcal{J}_0$ for every positive root $\alpha$.
\end{obs}

In \cite{AC1} the authors classified all flag manifold with two isotropy summands, this classification is presented in the following table:

\begin{table}[htb!]
\begin{center}
\begin{tabular}{|l|c|}
\hline $\mathbb{F}_\Theta = U/K$  & $\Sigma \backslash \Theta$ \\ \hline
$SO(2l+1)/U(l-m)\times SO(2l+1)$ \ ($l-m\not= 1$) & $\{ \alpha_{l-m} \}$ \\
$Sp(l)/ U(l-m)\times Sp(m)$ \ ($m\not= 0$) &  $\{ \alpha_{l-m} \}$ \\
$SO(2l)/ U(l-m)\times SO(2m)$ \ ($l-m\not=1$, $m\not= 0$) & $\{ \alpha_{l-m} \}$ \\
$G_2/U(2)$ \ ($U(2)$ represented by the short root) &  $\{ \alpha_{1} \}$ \\
$F_4/SO(7)\times U(1)$ &  $ \{ \alpha_{4} \}$\\
$F_4/Sp(3)\times U(1)$ &  $\{ \alpha_{1}\}$\\
$E_6/SU(6)\times U(1)$ &  $\{ \alpha_{6}\}$ \\
$E_6/SU(2)\times SU(5) \times U(1)$ &  $\{ \alpha_{2}\}$\\
$E_7/SU(7)\times U(1)$ &  $\{ \alpha_{7} \}$ \\
$E_7/SU(2)\times SO(10)\times U(1)$  & $\{ \alpha_{2}\}$ \\
$E_7/SO(12)\times U(1)$ & $\{ \alpha_{6}\}$ \\
$E_8/E_7\times U(1)$ &  $\{ \alpha_{1}\}$ \\
$E_8/SO(14)\times U(1)$ &  $\{ \alpha_{7}\}$ \\ \hline
\end{tabular}
\caption{Flag manifolds with two isotropy summands.}
\end{center}

\end{table}

Therefore, let $\mathcal{J}$ be an invariant generalized almost complex structure on $\mathbb{F}_\Theta$, where $\mathbb{F}_\Theta$ is one of the flag manifolds in the table above. Then $\mathcal{J}$ is integrable if and only if $\mathcal{J}_\alpha$ is of the same type for all $\alpha \in \Pi \backslash \langle \Theta \rangle$.
One example of this is presented now:\\

\begin{exem}
Let $\mathbb{F}_\Theta$ be a partial flag manifold of type $B_2$ generated by $\Theta = \{ \alpha \}$, where $\alpha \in \Sigma$ is the long root. In this case, 
\[
\mathfrak{m} = \mathfrak{m}_1 \oplus \mathfrak{m}_2,
\]
where $\mathfrak{m}_1 = \mathfrak{u}_\beta \oplus \mathfrak{u}_{\alpha+ \beta}$ and $\mathfrak{m}_2 = \mathfrak{u}_{\alpha + 2\beta}$. Let $\mathcal{J}$ be an invariant generalized almost complex structure on $\mathbb{F}_\Theta$. If $\mathcal{J}$ is integrable, we have that $\mathcal{J}_\beta$, $\mathcal{J}_{\alpha+\beta}$ and $\mathcal{J}_{\alpha+2\beta}$ are of the same type, complex or noncomplex. 
\end{exem} 

The second case is when $\alpha \in \mathfrak{m}_i$, $\beta \in \mathfrak{m}_j$ and $\alpha+\beta \in \mathfrak{m}_k$, with $i,j,k$ mutually distinct. In this case we have more possibilities, actually we have all the possible rows from Table \ref{Jintegravel}.

Our aim is to classify all invariant generalized complex structures on flag manifold with two, three and four isotropy summands. The case with two isotropy summands we have already done. In the papers \cite{Kim} and \cite{AC2} there is a complete classification of the partial flag manifolds with three and four isotropy summands, respectively. We present this classification on the following table: 

\begin{table}[htb!]
\begin{center}
\begin{tabular}{|l|c|c|}
\hline $\mathbb{F}_\Theta = U/K$  & $\mathfrak{m} = \oplus_{i=1} ^s \mathfrak{m}_i$ & $\Sigma \backslash \Theta$ \\ \hline
$SU(l_1 + l_2 + l_3)/S(U(l_1)\times U(l_2)\times U(l_3))$ & $s=3$ & $\{ \alpha_{l_1},\alpha_{l_2}\}$\\
$SO(2l)/U(1)\times U(l-1)$ \ ($l\geq 4$) & $s=3$ & $\{ \alpha_{l-1},\alpha_l\}$\\
$G_2/U(2)$ \ ($U(2)$ represented by the long root) & $s=3$ & $\{ \alpha_{2}\}$ \\
$F_4/U(2)\times SU(3)$ & $s=3$ & $\{ \alpha_{2}\}$\\
$E_6/U(1)\times U(1) \times SO(8)$ & $s=3$ & $\{ \alpha_{1},\alpha_5\}$\\
$E_6/U(2)\times SU(3) \times SU(3)$ & $s=3$ & $\{ \alpha_{3}\}$\\
$E_7/U(3)\times SU(5)$ & $s=3$ & $\{ \alpha_{3}\}$\\
$E_7/U(2)\times SU(6)$ & $s=3$ & $\{ \alpha_{5}\}$\\
$E_8/U(2)\times E_6$ & $s=3$ & $\{ \alpha_{2}\}$\\
$E_8/U_8$ & $s=3$ & $\{ \alpha_{8}\}$ \\ \hline
$F_4/SU(3)\times SU(2) \times SU(1)$ & $s=4$ & $\{ \alpha_{3}\}$ \\
$E_7/ SU(4)\times SU(3)\times SU(2) \times SU(1)$ & $s=4$ & $\{ \alpha_{4}\}$\\
$E_8/ SU(7)\times SU(2) \times U(1)$ & $s=4$ & $\{ \alpha_{6}\}$\\
$E_8/ SO(10) \times SU(3)\times U(1)$ & $s=4$ & $\{ \alpha_{3}\}$\\
$E_6/ SU(5) \times U(1)\times U(1)$ & $s=4$ & $\{ \alpha_{1},\alpha_2\}$\\
$E_7/ SO(10)\times U(1)\times U(1)$ &$s=4$ & $\{ \alpha_{1},\alpha_2\}$\\
$SO(2l+1)/ U(1)\times U(1)\times SO(2l-3)$ \ ($l\geq 2$) & $s=4$ & $\{ \alpha_{1},\alpha_2\}$\\
$SO(2l)/U(1)\times U(1)\times SO(2l-4)$ \ ($l\geq 3$) & $s=4$ & $\{ \alpha_{1},\alpha_2\}$\\
$SO(2l)/U(p)\times U(l-p)$ \ ($l\geq 4$, $2\leq p\leq l-2$) & $s=4$ & $\{ \alpha_{p},\alpha_l\}$\\
$Sp(l)/ U(p)\times U(l-p)$ \ ($l\geq 2$, $1\leq p\leq l-1$) & $s=4$ & $\{ \alpha_{p},\alpha_l\}$ \\ \hline
\end{tabular} 
\caption{Partial flag manifolds with three and four isotropy summands.}
\label{tabelaflag}
\end{center}
\end{table}

\noindent where the second column of the Table \ref{tabelaflag} indicates the number of irreducible components of the isotropy representation and the third column give the set $\Sigma \backslash \Theta$, where $\Theta$ is the subset of $\Sigma$ which defines $\mathbb{F}_\Theta$.

Recall that the height of a simple root $\alpha_i$ is the positive integer $m_i$ that appears on the highest root $\tilde{\alpha} = \sum_{j=1} ^l m_j\alpha_j$ of $\Pi^+$. Define the function $\hht\colon \Sigma \rightarrow \mathbb{Z}$, $\hht(\alpha_i) = m_i$. In \cite{AC1} the authors classified all flag manifolds with two isotropy summands. This was done considering $\Sigma \backslash \Theta = \{ \alpha_i \ | \ \hht(\alpha_i) = 2\}$. In \cite{Kim}, the author obtained all flag manifold with three isotropy summands, by setting $\Sigma \backslash \Theta = \{ \alpha_i \ | \ \hht(\alpha_i) = 3\}$ or $\Sigma \backslash \Theta = \{ \alpha_i, \alpha_j \ | \ \hht(\alpha_i) = \hht(\alpha_j) = 1\}$. The classification of all flag manifolds with four isotropy summands given in \cite{AC2} was done considering $\Sigma \backslash \Theta = \{ \alpha_i \ | \ \hht(\alpha_i) = 4\}$ or $\{ \alpha_i,\alpha_j \ | \ \hht(\alpha_i) = 1 \ \textrm{ and } \ \hht(\alpha_j)=2 \}$. \\

\begin{obs}
Corollary \ref{2somandos} characterize all invariant generalized complex structures on flag manifolds with two isotropy summands. 
\end{obs}

If the flag manifold $\mathbb{F}_\Theta$ is such that $\Sigma \backslash \Theta = \{ \alpha_i\}$ we can prove a result analogous to Corollary \ref{2somandos}. But first, let us prove the following:

\begin{lema}\label{ncomponents}
Let $\mathbb{F}_\Theta$ be a partial flag manifold where $\Theta = \Sigma \backslash \{ \alpha_{i_0}\}$. Suppose that $\hht(\alpha_{i_0}) = n$, then $\mathfrak{m}$ is decomposed into $n$ mutually non equivalent irreducible components, that is, 
\[
\mathfrak{m} = \mathfrak{m}_1 \oplus \cdots \oplus \mathfrak{m}_n.
\]
\end{lema}
\begin{demons}
Let $\mathbb{F}_\Theta$ be a flag manifold where $\Sigma \backslash \Theta = \{ \alpha_{i_0}\}$ for some $i_0 \in \{ 1,\cdots, l\}$. Let $H_0 = a_{i_0}H_{i_0}$ be a generic element of the positive Weyl chamber. Since $\hht(\alpha_{i_0}) = n$ we will have exactly $n$ distinct eigenvalues of $\ad (H_0)$. In fact, given $j\in \{1,\cdots ,n\}$ we always have at least one root $\alpha \in \Pi^+$ such that $\alpha = \sum_{k=1} ^l n_k \alpha_k$ with $n_{i_0} =j$ and therefore $\alpha (H_0) = ja_{i_0}$, proving that $\ad (H_0)$ has $n$ eigenvalues. Since the number of irreducible components is the number of eigenvalues of $\ad (H_0)$, we have proved that $\mathfrak{m}$ is decomposed into $n$ irreducible components.   
\end{demons}

Using this result we can prove the following.

\begin{prop}\label{oneroot}
Let $\mathbb{F}_\Theta$ be a partial flag manifold where $\Theta = \Sigma \backslash \{ \alpha_{i_0}\}$. Consider $\mathcal{J}$ an invariant generalized almost complex structure on $\mathbb{F}_\Theta$. Then $\mathcal{J}$ is integrable if and only if $\mathcal{J}_\alpha$ is of the same type for all $\alpha \in \Pi \backslash \langle \Theta \rangle$.
\end{prop}
\begin{demons}
Let $\mathbb{F}_\Theta$ be a partial flag manifold where $\Theta = \Sigma \backslash \{ \alpha_{i_0}\}$ and suppose that $\hht (\alpha_{i_0}) = n \geq 3$, because $n=1$ is trivial and $n=2$ is given by Corollary \ref{2somandos}. Then, by Lemma \ref{ncomponents}, we have that $\mathfrak{m} = \mathfrak{m}_1 \oplus \cdots \oplus \mathfrak{m}_n$. Set $\mathfrak{m}_j$ the component associated to the eigenvalue $\alpha (H_0) = ja_{i_0}$, where $H_0 = a_{i_0}H_{i_0}$ is a generic element of the positive Weyl chamber. Thus given $\alpha,\beta \in \mathfrak{m}_1$ such that $\alpha+\beta$ is a root, then we must have $\alpha+\beta \in \mathfrak{m}_2$, because $(\alpha + \beta)(H_0) = 2a_{i_0}$. Thus, if $\mathcal{J}$ is an invariant generalized (integrable) complex structure on $\mathbb{F}_\Theta$, by Table \ref{Jintegravel}, we must have $\mathcal{J}|_{\mathfrak{m}_1 \oplus \mathfrak{m}_1 ^\ast}$ and $\mathcal{J}|_{\mathfrak{m}_2 \oplus \mathfrak{m}_2 ^\ast}$ of same type. Suppose that for $j\geq 1$ we have $\mathcal{J}|_{\mathfrak{m}_1 \oplus \mathfrak{m}_1 ^\ast},\cdots ,\mathcal{J}|_{\mathfrak{m}_j \oplus \mathfrak{m}_j ^\ast}$ are of the same type and let us prove that $\mathcal{J}|_{\mathfrak{m}_{j+1} \oplus \mathfrak{m}_{j+1} ^\ast}$ must be of the same type. Let $\alpha \in \mathfrak{m}_1$ and $\beta \in \mathfrak{m}_j$ such that $\alpha+\beta$ is also a root, then we must have that $\alpha+\beta \in \mathfrak{m}_{j+1}$. Since $\mathcal{J}|_{\mathfrak{m}_1 \oplus \mathfrak{m}_1 ^\ast}$ and $\mathcal{J}|_{\mathfrak{m}_j \oplus \mathfrak{m}_j ^\ast}$ have the same type, follows from Table \ref{Jintegravel} that $\mathcal{J}|_{\mathfrak{m}_{j+1} \oplus \mathfrak{m}_{j+1} ^\ast}$ has the same type of them.

\end{demons}


Observe that Proposition \ref{oneroot} covers most cases in Table \ref{tabelaflag}, it remains only nine cases to analyse. We will do it case by case, analysing each type of algebra associated to $\mathbb{F}_\Theta$.

\subsubsection*{Type $A_l$}  
The Lie algebra of type $A_l$ has Dynkin diagram
\begin{center}
  \begin{tikzpicture}[scale=.4]
    \draw (0,0) circle (.3cm) node [below] {$\alpha_1$};
    \draw (0.3 cm,0) -- +(1.4 cm,0); 
    \draw (2,0) circle (.3cm) node [below] {$\alpha_2$};
    \draw (2.3 cm,0) -- +(1.4 cm,0); 
    \draw[dotted] (3.7cm,0) -- +(3cm,0);
    \draw (6.3 cm,0) -- +(1.4 cm,0); 
    \draw (8,0) circle (.3cm) node [below] {$\alpha_{l-1}$};
    \draw (8.3 cm,0) -- +(1.4 cm,0); 
    \draw (10,0) circle (.3cm) node [below] {$\alpha_{l}$};
  \end{tikzpicture}
\end{center}
and it is associated to the algebra $\mathfrak{sl}(l+1)$. A Cartan subalgebra is the diagonal matrices with trace zero. The roots are $\lambda_i - \lambda_j$, with $i\not= j$, where $\lambda_i$ is given by 
\[
\lambda_i\colon \textnormal{diag}\{a_1, \cdots ,a_{l+1}\} \longmapsto a_i.
\]

A simple root system is given by
\[
\Sigma = \{ \lambda_1 - \lambda_2, \lambda_2 - \lambda_3, \cdots ,\lambda_{l} - \lambda_{l+1}\},
\]
and the positive roots are
\[
\Pi^+=\{ \lambda_i - \lambda_j \ | \ i<j\}.
\]
Using the notation $\Sigma = \{ \alpha_1, \alpha_2, \cdots, \alpha_l\}$, where $\alpha_i = \lambda_i - \lambda_{i+1}$, we can write the positive roots as linear combination of the simple roots by
\[
\{ \alpha_i + \alpha_{i+1} + \cdots + \alpha_j \ | \ 1 \leq i < j \leq l\}.
\]

In Table \ref{tabelaflag}, $A_l$ appears only in one case, to know, $\mathbb{F}_\Theta = SU(l_1+l_2+l_3)/S(U(l_1)\times U(l_2)\times U(l_3))$. In this case, we have $\Theta = \Sigma \backslash \{\alpha_{l_1} , \alpha_{l_2}\}$ in which $\mathfrak{m}$ has three irreducible components $\mathfrak{m}_1$, $\mathfrak{m}_2$ and $\mathfrak{m}_3$. Thus the components are described by:
\begin{eqnarray*}
\mathfrak{m}_1 & = & \{ \alpha_i +\cdots + \alpha_{l_1} + \cdots + \alpha_j \ | \ 1\leq i \leq l_1 \leq j \leq l_2 -1\}
\end{eqnarray*}
\begin{eqnarray*}
\mathfrak{m}_2 & = & \{ \alpha_i + \cdots +\alpha_{l_2} + \cdots + \alpha_j \ | \ l_1 + 1 \leq i \leq l_2 \leq j \leq l\}
\end{eqnarray*}
\begin{eqnarray*}
\mathfrak{m}_3 & = & \{ \alpha_i + \cdots + \alpha_{l_1} + \cdots + \alpha_{l_2} + \cdots + \alpha_j \ | \ 1 \leq i \leq l_1 \ \textrm{ e } \ l_2 \leq j \leq l\}
\end{eqnarray*}
where $l = l_1+l_2+l_3$. Therefore, given a triple of roots ($\alpha, \beta, \alpha+\beta$), it is immediate to see that $\alpha \in \mathfrak{m}_1$, $\beta \in \mathfrak{m}_2$ and $\alpha+\beta \in \mathfrak{m}_3$. Thus, given $\mathcal{J}$ an invariant generalized almost complex structure on $\mathbb{F}_\Theta$, we have that if $\mathcal{J}$ is integrable then $\mathcal{J}$ is given by Table \ref{Jintegravel}, that is, is one of the rows on the following table:
\begin{center}
\begin{tabular}{c|c|c}
$\mathcal{J}|_{\mathfrak{m}_1 \oplus \mathfrak{m}_1 ^\ast}$ & $\mathcal{J}|_{\mathfrak{m}_2 \oplus \mathfrak{m}_2 ^\ast}$ & $\mathcal{J}|_{\mathfrak{m}_3 \oplus \mathfrak{m}_3 ^\ast}$ \\ \hline
complex ($\pm \mathcal{J}_0$) & complex ($\pm \mathcal{J}_0$) & complex ($\pm \mathcal{J}_0$) \\
complex ($\pm \mathcal{J}_0$) & complex ($ \mp \mathcal{J}_0$) & complex ($\pm \mathcal{J}_0$) \\
complex ($\pm \mathcal{J}_0$) & complex ($\mp \mathcal{J}_0$) & complex ($\mp \mathcal{J}_0$) \\
complex ($\pm \mathcal{J}_0$) & noncomplex & complex ($\pm \mathcal{J}_0$) \\
noncomplex  & complex ($\pm \mathcal{J}_0$) & complex ($\pm \mathcal{J}_0$) \\
complex ($\pm \mathcal{J}_0$) & complex ($\mp \mathcal{J}_0$) & noncomplex \\
noncomplex & noncomplex & noncomplex
\end{tabular} 
\end{center}

\subsubsection*{Type $B_l$}
The diagram
\begin{center}
  \begin{tikzpicture}[scale=.4]
    \draw (-1,0) node[anchor=east]  {$B_l, l\geq 2$};
    \draw (0,0) circle (.3cm) node [below] {$\alpha_1$};
    \draw (0.3 cm,0) -- +(1.4 cm,0); 
    \draw (2,0) circle (.3cm) node [below] {$\alpha_2$};
    \draw (2.3 cm,0) -- +(1.4 cm,0); 
    \draw[dotted] (3.7cm,0) -- +(3cm,0);
    \draw (6.3 cm,0) -- +(1.4 cm,0); 
    \draw (8,0) circle (.3cm) node [below] {$\alpha_{l-1}$};
    \draw (8.3 cm, .1 cm) -- +(1.4 cm,0);
    \draw (8.3 cm, -.1 cm) -- +(1.4 cm,0);
    \draw (10,0) circle (.3cm) node [below] {$\alpha_{l}$};
    \draw (9.7 cm,0cm) -- (9.3 cm,.8cm);
    \draw (9.7 cm,0cm) -- (9.3 cm,-0.8cm);
  \end{tikzpicture}
\end{center}
is the Dynkin diagram of the algebra $\mathfrak{so}(2l+1) = \{ A \in \mathfrak{sl}(2l+1) \ | \ A+A^t=0\}$
whose dimension is $\frac{(2l+1)2l}{2}=l(2l+1)$. A Cartan subalgebra $\mathfrak{h}$ is the one dimensional subalgebra of the diagonal matrices in $\mathfrak{so}(2l+1)$, that is, $H\in \mathfrak{h}$ if and only  if $H$ is of the form
\[
H = \left( \begin{array}{ccc}
0 & & \\
& \Lambda & \\
& & -\Lambda
\end{array}\right)
\]
with $\Lambda$ being an arbitrary diagonal $l\times l$ matrix. Let $\lambda_i$, $i=1,\cdots ,l$, be the functional
\[
\lambda_i\colon \textnormal{diag}(a_1, \cdots , a_l)\longmapsto a_i.
\]
The roots are given by:
\begin{itemize}
\item $\pm \lambda_j$ with $j=1,\cdots, l$;

\item $\pm (\lambda_i - \lambda_j)$ with $1\leq i<j \leq l$;

\item $\pm (\lambda_i + \lambda_j)$ with $1 \leq i\not= j \leq l$. 
\end{itemize}
A simple root system is given by 
\[
\Sigma = \{ \lambda_1 - \lambda_2, \cdots, \lambda_{l-1} - \lambda_l, \lambda_l \}.
\]
The positive roots are $\lambda_j$ for $j=1,\cdots, l$, $\lambda_i - \lambda_j$ with $i<j$ and $\lambda_i + \lambda_j$ with $i\not=j$. If we use the notation $\Sigma =\{ \alpha_1, \cdots, \alpha_{l-1}, \alpha_l\}$, the positive roots can be written as
\begin{eqnarray*}
\lambda_j & = & \alpha_j + \alpha_{j+1} + \cdots + \alpha_l \\
\lambda_i+\lambda_j & = & \alpha_i + \cdots + \alpha_{j-1} + 2\alpha_j + \cdots + 2 \alpha_l \\
\lambda_i-\lambda_j & = & \alpha_i + \cdots + \alpha_{j-1}.
\end{eqnarray*}

The type $B_l$ also appears only once in Table $\ref{tabelaflag}$, to know, $\mathbb{F}_\Theta = SO(2l+1)/ U(1)\times U(1) \times SO(2l-3)$ with $l\geq 2$, where $\Theta = \Sigma \backslash \{ \alpha_1,\alpha_2\}$. In this case, $\mathfrak{m}$ has four irreducible components, that is, $\mathfrak{m} = \mathfrak{m}_1 \oplus \mathfrak{m}_2 \oplus \mathfrak{m}_3 \oplus \mathfrak{m}_4$ in which
\begin{eqnarray*}
\mathfrak{m}_1 & = & \{ \alpha_1\} 
\end{eqnarray*}
\begin{eqnarray*}
\mathfrak{m}_2 & = & \{ \alpha_2 + \cdots + \alpha_i , \alpha_2 + \cdots + \alpha_{i-1} + 2\alpha_i + \cdots + 2 \alpha_l \ | \ 2 \leq i\leq l \}
\end{eqnarray*}
\begin{eqnarray*}
\mathfrak{m}_3 & = & \{ \alpha_1 + \alpha_2 + \cdots + \alpha_i , \alpha_1 +\alpha_2 + \cdots + \alpha_{i-1} + 2\alpha_i + \cdots + 2\alpha_l\ | \ 3 \leq i \leq l \}
\end{eqnarray*}
\begin{eqnarray*}
\mathfrak{m}_4 & = & \{ \alpha_1 + 2\alpha_2 + \cdots + 2\alpha_2 \}. 
\end{eqnarray*}

Now, observe that if $\alpha_1 \in \mathfrak{m}_1$ and $\beta \in \mathfrak{m}_2$ then $\alpha_1 + \beta \in \mathfrak{m}_3$. Analogously, it is possible to see that there are roots $\beta \in \mathfrak{m}_2$ and $\gamma \in \mathfrak{m}_3$ such that $\beta +\gamma \in \mathfrak{m}_4$. Thus, given $\mathcal{J}$ an invariant generalized almost complex structure on $\mathbb{F}_\Theta$, we have that if $\mathcal{J}$ is integrable, then it is one of the following possibilities:
\begin{center}
\begin{tabular}{c|c|c|c}
$\mathcal{J}|_{\mathfrak{m}_1 \oplus \mathfrak{m}_1 ^\ast}$ & $\mathcal{J}|_{\mathfrak{m}_2 \oplus \mathfrak{m}_2 ^\ast}$ & $\mathcal{J}|_{\mathfrak{m}_3 \oplus \mathfrak{m}_3 ^\ast}$ & $\mathcal{J}|_{\mathfrak{m}_4 \oplus \mathfrak{m}_4 ^\ast}$  \\ \hline
complex ($\pm \mathcal{J}_0$) & complex ($\pm \mathcal{J}_0$) & complex ($\pm \mathcal{J}_0$) & complex ($\pm \mathcal{J}_0$) \\
complex ($\pm \mathcal{J}_0$) & complex ($\mp \mathcal{J}_0$) & complex ($\pm \mathcal{J}_0$) & complex ($\pm \mathcal{J}_0$) \\
complex ($\pm \mathcal{J}_0$) & complex ($\mp \mathcal{J}_0$) & complex ($\pm \mathcal{J}_0$) & complex ($\mp \mathcal{J}_0$) \\
complex ($\pm \mathcal{J}_0$) & complex ($\mp \mathcal{J}_0$) & complex ($\mp \mathcal{J}_0$) & complex ($\mp \mathcal{J}_0$) \\
complex ($\pm \mathcal{J}_0$) & complex ($\mp \mathcal{J}_0$) & complex ($\pm \mathcal{J}_0$) & noncomplex \\
complex ($\pm \mathcal{J}_0$) & complex ($\mp \mathcal{J}_0$) & noncomplex  & complex ($\mp \mathcal{J}_0$) \\
complex ($\pm \mathcal{J}_0$) & noncomplex  & complex ($\pm \mathcal{J}_0$) & complex ($\pm \mathcal{J}_0$) \\
noncomplex & complex ($\pm \mathcal{J}_0$) & complex ($\pm \mathcal{J}_0$) & complex ($\pm \mathcal{J}_0$) \\
noncomplex & noncomplex & noncomplex & noncomplex
\end{tabular} 
\end{center}

\subsubsection*{Type $C_l$}
The diagram
\begin{center}
  \begin{tikzpicture}[scale=.4]
    \draw (-1,0) node[anchor=east]  {$C_l, l\geq 3$};
    \draw (0,0) circle (.3cm) node [below] {$\alpha_1$};
    \draw (0.3 cm,0) -- +(1.4 cm,0); 
    \draw (2,0) circle (.3cm) node [below] {$\alpha_2$};
    \draw (2.3 cm,0) -- +(1.4 cm,0); 
    \draw[dotted] (3.7cm,0) -- +(3cm,0);
    \draw (6.3 cm,0) -- +(1.4 cm,0); 
    \draw (8,0) circle (.3cm) node [below] {$\alpha_{l-1}$};
    \draw (8.3 cm, .1 cm) -- +(1.4 cm,0);
    \draw (8.3 cm, -.1 cm) -- +(1.4 cm,0);
    \draw (10,0) circle (.3cm) node [below] {$\alpha_{l}$};
    \draw (8.3 cm,0cm) -- (8.7 cm,.8cm);
    \draw (8.3 cm,0cm) -- (8.7 cm,-0.8cm);
  \end{tikzpicture}
\end{center}
is the Dynkin diagram associated to $\mathfrak{sp}(l) = \{ A\in \mathfrak{sl}(2l) \ | \ AJ + JA^t = 0\}$, where 
\[
J = \left( \begin{array}{cc}
0 & -1 \\
1 & 0
\end{array} \right).
\] 
A Cartan subalgebra $\mathfrak{h}$ is the subalgebra of diagonal matrices in $\mathfrak{sp}(l)$, that is, the elements $H\in \mathfrak{h}$ are of the form 
\[
H = \left( \begin{array}{cc}
\Lambda & 0\\
0 & -\Lambda
\end{array} \right)
\]
with $\Lambda = \textnormal{diag}(a_1,\cdots ,a_l)$ an arbitrary diagonal $l\times l$ matrix. Thus, the roots are given by:
\begin{itemize}
\item $\pm (\lambda_i - \lambda_j)$ with $i < j$;

\item $\pm (\lambda_i + \lambda_j)$ for $i,j=1,\cdots ,l$. In particular, we have that $\pm 2\lambda_i$ is a root.
\end{itemize}
A simple root system is given by 
\[
\Sigma = \{ \lambda_1-\lambda_2,\cdots ,\lambda_{l-1}-\lambda_l,2\lambda_l\}.
\]
The positive roots are $\lambda_i - \lambda_j$ with $i<j$ e $\lambda_i +\lambda_j$. Using the notation $\Sigma = \{ \alpha_1, \cdots ,\alpha_{l-1},\alpha_l \}$, we have that the positive roots can be written as:
\begin{eqnarray*}
\lambda_i - \lambda_j & = & \alpha_i +\cdots + \alpha_{j-1}\\
\lambda_i + \lambda_j & = & \alpha_i + \cdots + \alpha_{j-1} + 2\alpha_j + \cdots + 2\alpha_{l-1} + \alpha_l.
\end{eqnarray*}

Analogous to the $B_l$ case, we have one flag manifold in Table \ref{tabelaflag}, this flag manifold is given by $\mathbb{F}_\Theta = Sp(l)/U(p)\times U(l-p)$ with $l\geq 2$ and $2\leq p \leq l-2$, where $\Theta = \Sigma \backslash \{ \alpha_p , \alpha_l\}$. In this case we have $\mathfrak{m} = \mathfrak{m}_1 \oplus \mathfrak{m}_2 \oplus \mathfrak{m}_3 \oplus \mathfrak{m}_4$. To know,
\begin{eqnarray*}
\mathfrak{m}_1 & = & \{ \alpha_l , \alpha_i +\cdots + \alpha_{j-1} + 2\alpha_j + \cdots + 2\alpha_{l-1} + \alpha_l \ | \ p < i \leq j < l\}
\end{eqnarray*}
\begin{eqnarray*}
\mathfrak{m}_2 & = & \{ \alpha_i + \cdots + \alpha_p + \cdots + \alpha_j \ | \ 1\leq i\leq j\leq l-1\}
\end{eqnarray*}
\begin{eqnarray*}
\mathfrak{m}_3 & = & \{ \alpha_i + \cdots + \alpha_p + \cdots + \alpha_{j-1} + 2\alpha_j + \cdots + 2\alpha_{l-1} + \alpha_j \ | \\
& & \ 1\leq i \leq p < j \leq l-1 \}
\end{eqnarray*}
\begin{eqnarray*}
\mathfrak{m}_4 & = & \{ \alpha_i + \cdots + \alpha_{j-1} + 2\alpha_j+ \cdots + 2\alpha_p + \cdots + 2\alpha_{l-1} + \alpha_l \ | \ 1\leq i \leq j \leq p\}
\end{eqnarray*}
It is easy to see that if $\alpha \in \mathfrak{m}_1$ and $\beta \in \mathfrak{m}_2$ with $\alpha+\beta$ root, then $\alpha + \beta \in \mathfrak{m}_3$. And, analogously, if $\alpha \in \mathfrak{m}_2$ and $\beta \in \mathfrak{m}_3$ with $\alpha+\beta$ root, then $\alpha + \beta \in \mathfrak{m}_4$. This way we obtain a case analogous to $B_l$, that is, if $\mathcal{J}$ is an invariant integrable generalized complex structure on $\mathbb{F}_\Theta$, then it must be one of the following cases:
\begin{center}
\begin{tabular}{c|c|c|c}
$\mathcal{J}|_{\mathfrak{m}_1 \oplus \mathfrak{m}_1 ^\ast}$ & $\mathcal{J}|_{\mathfrak{m}_2 \oplus \mathfrak{m}_2 ^\ast}$ & $\mathcal{J}|_{\mathfrak{m}_3 \oplus \mathfrak{m}_3 ^\ast}$ & $\mathcal{J}|_{\mathfrak{m}_4 \oplus \mathfrak{m}_4 ^\ast}$  \\ \hline
complex ($\pm \mathcal{J}_0$) & complex ($\pm \mathcal{J}_0$) & complex ($\pm \mathcal{J}_0$) & complex ($\pm \mathcal{J}_0$) \\
complex ($\pm \mathcal{J}_0$) & complex ($\mp \mathcal{J}_0$) & complex ($\pm \mathcal{J}_0$) & complex ($\pm \mathcal{J}_0$) \\
complex ($\pm \mathcal{J}_0$) & complex ($\mp \mathcal{J}_0$) & complex ($\pm \mathcal{J}_0$) & complex ($\mp \mathcal{J}_0$) \\
complex ($\pm \mathcal{J}_0$) & complex ($\mp \mathcal{J}_0$) & complex ($\mp \mathcal{J}_0$) & complex ($\mp \mathcal{J}_0$) \\
complex ($\pm \mathcal{J}_0$) & complex ($\mp \mathcal{J}_0$) & complex ($\pm \mathcal{J}_0$) & noncomplex \\
complex ($\pm \mathcal{J}_0$) & complex ($\mp \mathcal{J}_0$) & noncomplex  & complex ($\mp \mathcal{J}_0$) \\
complex ($\pm \mathcal{J}_0$) & noncomplex  & complex ($\pm \mathcal{J}_0$) & complex ($\pm \mathcal{J}_0$) \\
noncomplex & complex ($\pm \mathcal{J}_0$) & complex ($\pm \mathcal{J}_0$) & complex ($\pm \mathcal{J}_0$) \\
noncomplex & noncomplex & noncomplex & noncomplex
\end{tabular} 
\end{center}

\subsubsection*{Type $D_l$}

The diagram
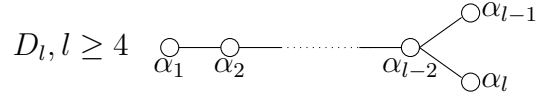
\begin{center}
  \begin{tikzpicture}[scale=.4]
    \draw (-1,0) node[anchor=east]  {$D_l, l\geq 4$};
    \draw (0,0) circle (.3cm) node [below] {$\alpha_1$};
    \draw (0.3 cm,0) -- +(1.4 cm,0); 
    \draw (2,0) circle (.3cm) node [below] {$\alpha_2$};
    \draw (2.3 cm,0) -- +(1.4 cm,0); 
    \draw[dotted] (3.7cm,0) -- +(3cm,0);
    \draw (6.3 cm,0) -- +(1.4 cm,0); 
    \draw (8,0) circle (.3cm) node [below] {$\alpha_{l-2}$};
    \draw (8.3 cm,0) -- +(1.4 cm,1cm);
    \draw (8.3 cm,0) -- +(1.4 cm,-1cm); 
    \draw (10cm,1.15cm) circle (.3cm) node[anchor=west] {$\alpha_{l-1}$};
    \draw (10cm,-1.15cm) circle (.3cm) node [anchor=west] {$\alpha_{l}$};
  \end{tikzpicture}
\end{center}
is the Dynkin diagram associated to $\mathfrak{so}(2l) = \{ A \in \mathfrak{sl}(2l) \ | \ A+A^t = 0\}$. A Cartan subalgebra $\mathfrak{h}$ is the subalgebra of diagonal matrices, where its elements are written as
\[
H = \left( \begin{array}{cc}
\Lambda & 0 \\
0 & -\Lambda
\end{array}\right)
\]
with $\Lambda = \textnormal{diag}(a_1,\cdots ,a_l)$ an arbitrary diagonal matrix. The roots are given by:
\begin{itemize}
\item $\pm (\lambda_i - \lambda_j)$ with $i < j$;
\item $\pm (\lambda_i + \lambda_j)$ with $i \not= j$.
\end{itemize}
A simple root system is given by:
\[
\Sigma = \{ \lambda_1 - \lambda_2,\cdots, \lambda_{l-1} - \lambda_l, \lambda_{l-1} + \lambda_l\}.
\]
The positive roots are $\lambda_i - \lambda_j$ com $i<j$ e $\lambda_i + \lambda_j$ with $i\not= j$. Using the notation $\Sigma = \{ \alpha_1, \cdots ,\alpha_{l-1}, \alpha_l \}$, the positive roots can be written as
\begin{eqnarray*}
\lambda_i - \lambda_j & = & \alpha_i + \cdots + \alpha_{j-1}\\
\lambda_i + \lambda_j & = & \alpha_i + \cdots + \alpha_{j-1} + 2\alpha_j + \cdots + 2\alpha_{l-2} + \alpha_{l-1} + \alpha_l.
\end{eqnarray*}

In Table \ref{tabelaflag} the type $D_l$ appears three times:
\begin{enumerate}
\item $SO(2l)/ U(1)\times U(l-1)$ with $l\geq 4$;

\item $SO(2l)/ U(1)\times U(1) \times SO(2l-4)$ with $l\geq 3$;

\item $SO(2l)/ U(p)\times U(l-p)$ with $l\geq 4$ and $2\leq p\leq l-2$.
\end{enumerate}

The first case is a flag manifold with three summands. In this case, we have $\mathbb{F}_\Theta$ in which $\Theta = \Sigma \backslash \{ \alpha_{l-1},\alpha_l \}$ and $\mathfrak{m} = \mathfrak{m}_1 \oplus \mathfrak{m}_2 \oplus \mathfrak{m}_3$ with
\begin{eqnarray*}
\mathfrak{m}_1 & = & \{ \alpha_i + \cdots + \alpha_l-1 \ | \ 1\leq i \leq l-1 \}
\end{eqnarray*}
\begin{eqnarray*}
\mathfrak{m}_2 & = & \{ \alpha_l, \alpha_i + \cdots + \alpha_{l-2} + \alpha_l \ | \ 1\leq i \leq l-2\}
\end{eqnarray*}
\begin{eqnarray*}
\mathfrak{m}_3 & = & \{ \alpha_1 + \cdots + \alpha_{l-1} + \alpha_l , \alpha_i + \cdots + \alpha_{j-1} + 2\alpha_j + \cdots + 2\alpha_{l-2} + \alpha_{l-1} + \alpha_l \ | \\ 
& & \ 1\leq i < j \leq l-2\}.
\end{eqnarray*}
Thus, it is easy to see that if $\alpha \in \mathfrak{m}_1$ and $\beta \in \mathfrak{m}_2$ such that $\alpha+\beta$ is a root, then $\alpha + \beta \in \mathfrak{m}_3$. Therefore, given an invariant generalized almost complex structure $\mathcal{J}$, it is integrable if and only if it is one of the following combinations:
\begin{center}
\begin{tabular}{c|c|c}
$\mathcal{J}|_{\mathfrak{m}_1 \oplus \mathfrak{m}_1 ^\ast}$ & $\mathcal{J}|_{\mathfrak{m}_2 \oplus \mathfrak{m}_2 ^\ast}$ & $\mathcal{J}|_{\mathfrak{m}_3 \oplus \mathfrak{m}_3 ^\ast}$ \\ \hline
complex ($\pm \mathcal{J}_0$) & complex ($\pm \mathcal{J}_0$) & complex ($\pm \mathcal{J}_0$) \\
complex ($\pm \mathcal{J}_0$) & complex ($ \mp \mathcal{J}_0$) & complex ($\pm \mathcal{J}_0$) \\
complex ($\pm \mathcal{J}_0$) & complex ($\mp \mathcal{J}_0$) & complex ($\mp \mathcal{J}_0$) \\
complex ($\pm \mathcal{J}_0$) & noncomplex & complex ($\pm \mathcal{J}_0$) \\
noncomplex  & complex ($\pm \mathcal{J}_0$) & complex ($\pm \mathcal{J}_0$) \\
complex ($\pm \mathcal{J}_0$) & complex ($\mp \mathcal{J}_0$) & noncomplex \\
noncomplex & noncomplex & noncomplex
\end{tabular} 
\end{center}

In the second case, we have $\mathbb{F}_\Theta$ where $\Theta = \Sigma \backslash \{\alpha_1 , \alpha_2\}$ and $\mathfrak{m} = \mathfrak{m}_1 \oplus \mathfrak{m}_2 \oplus \mathfrak{m}_3 \oplus \mathfrak{m}_4$ with
\begin{eqnarray*}
\mathfrak{m}_1 & = & \{ \alpha_1 \}
\end{eqnarray*}
\begin{eqnarray*}
\mathfrak{m}_2 & = & \{ \alpha_2 + \cdots + \alpha_i \ | \ 3\leq i \leq l\} \cup \\ 
& & \{ \alpha_2 + \cdots + \alpha_{j-1} + 2\alpha_j + \cdots + 2\alpha_{l-2}+\alpha_{l-1} + \alpha_j \ | \ 2<j\leq l-2\}
\end{eqnarray*}
\begin{eqnarray*}
\mathfrak{m}_3 & = & \{ \alpha_1 + \alpha_2 + \cdots + \alpha_i \ | \ i\geq 3\} \cup \\ 
& & \{ \alpha_1 + \alpha_2 + \cdots + \alpha_{j-1} + 2\alpha_j + \cdots + 2\alpha_{l-2}+\alpha_{l-1} + \alpha_j \ | \ 2<j\leq l-2\}
\end{eqnarray*}
\begin{eqnarray*}
\mathfrak{m}_4 & = & \{ \alpha_1 + 2\alpha_2 + \cdots + 2\alpha_{l-2} + \alpha_{l-1} + \alpha_l \}.
\end{eqnarray*}
Thus, if $\mathcal{J}$ is an invariant integrable generalized complex structure, it must be one of the following structures:
\begin{center}
\begin{tabular}{c|c|c|c}
$\mathcal{J}|_{\mathfrak{m}_1 \oplus \mathfrak{m}_1 ^\ast}$ & $\mathcal{J}|_{\mathfrak{m}_2 \oplus \mathfrak{m}_2 ^\ast}$ & $\mathcal{J}|_{\mathfrak{m}_3 \oplus \mathfrak{m}_3 ^\ast}$ & $\mathcal{J}|_{\mathfrak{m}_4 \oplus \mathfrak{m}_4 ^\ast}$  \\ \hline
complex ($\pm \mathcal{J}_0$) & complex ($\pm \mathcal{J}_0$) & complex ($\pm \mathcal{J}_0$) & complex ($\pm \mathcal{J}_0$) \\
complex ($\pm \mathcal{J}_0$) & complex ($\mp \mathcal{J}_0$) & complex ($\pm \mathcal{J}_0$) & complex ($\pm \mathcal{J}_0$) \\
complex ($\pm \mathcal{J}_0$) & complex ($\mp \mathcal{J}_0$) & complex ($\pm \mathcal{J}_0$) & complex ($\mp \mathcal{J}_0$) \\
complex ($\pm \mathcal{J}_0$) & complex ($\mp \mathcal{J}_0$) & complex ($\mp \mathcal{J}_0$) & complex ($\mp \mathcal{J}_0$) \\
complex ($\pm \mathcal{J}_0$) & complex ($\mp \mathcal{J}_0$) & complex ($\pm \mathcal{J}_0$) & noncomplex \\
complex ($\pm \mathcal{J}_0$) & complex ($\mp \mathcal{J}_0$) & noncomplex  & complex ($\mp \mathcal{J}_0$) \\
complex ($\pm \mathcal{J}_0$) & noncomplex  & complex ($\pm \mathcal{J}_0$) & complex ($\pm \mathcal{J}_0$) \\
noncomplex & complex ($\pm \mathcal{J}_0$) & complex ($\pm \mathcal{J}_0$) & complex ($\pm \mathcal{J}_0$) \\
noncomplex & noncomplex & noncomplex & noncomplex
\end{tabular} 
\end{center}

The last case is given by $\mathbb{F}_\Theta$ with $\Theta = \Sigma \backslash \{ \alpha_p, \alpha_l\}$, remembering that $2\leq p \leq l-2$. In this case $\mathfrak{m} = \mathfrak{m}_1 \oplus \mathfrak{m}_2 \oplus \mathfrak{m}_3 \oplus \mathfrak{m}_4$ with
\begin{eqnarray*}
\mathfrak{m}_1 & = & \{ \alpha_j +\cdots + \alpha_l \ | \ p<j\leq l\} \cup \\
& & \{ \alpha_i + \cdots + \alpha_{j-1} +2\alpha_j + \cdots + 2\alpha_{l-2} + \alpha_{l-1} + \alpha_l \ | \ p<i\leq j \leq l-2\}
\end{eqnarray*}
\begin{eqnarray*}
\mathfrak{m}_2 & = & \{ \alpha_i + \cdots + \alpha_p + \cdots + \alpha_j \ | \ 1\leq i \leq p \leq j \leq l-1\}
\end{eqnarray*}
\begin{eqnarray*}
\mathfrak{m}_3 & = & \{\alpha_i + \cdots + \alpha_p + \cdots + \alpha_l \ | \ 1\leq i\leq p\} \cup \\ 
& & \{\alpha_i + \cdots + \alpha_p + \cdots + \alpha_{j-1} +2\alpha_j + \cdots + 2\alpha_{l-2} + \alpha_{l-1} + \alpha_l \ | \\
& & \ i\leq p \leq j \leq l-2\}
\end{eqnarray*}
\begin{eqnarray*}
\mathfrak{m}_4 & = & \{\alpha_i + \cdots + \alpha_{j-1} +2\alpha_j + \cdots + 2\alpha_p \cdots + 2\alpha_{l-2} + \alpha_{l-1} + \alpha_l \ | \\
& & \ i\leq p \leq j \leq l-2\}
\end{eqnarray*}
Therefore $\mathcal{J}$ is an invariant integrable generalized complex structure on $\mathbb{F}_\Theta$ if and only if it is one of the following cases:
\begin{center}
\begin{tabular}{c|c|c|c}
$\mathcal{J}|_{\mathfrak{m}_1 \oplus \mathfrak{m}_1 ^\ast}$ & $\mathcal{J}|_{\mathfrak{m}_2 \oplus \mathfrak{m}_2 ^\ast}$ & $\mathcal{J}|_{\mathfrak{m}_3 \oplus \mathfrak{m}_3 ^\ast}$ & $\mathcal{J}|_{\mathfrak{m}_4 \oplus \mathfrak{m}_4 ^\ast}$  \\ \hline
complex ($\pm \mathcal{J}_0$) & complex ($\pm \mathcal{J}_0$) & complex ($\pm \mathcal{J}_0$) & complex ($\pm \mathcal{J}_0$) \\
complex ($\pm \mathcal{J}_0$) & complex ($\mp \mathcal{J}_0$) & complex ($\pm \mathcal{J}_0$) & complex ($\pm \mathcal{J}_0$) \\
complex ($\pm \mathcal{J}_0$) & complex ($\mp \mathcal{J}_0$) & complex ($\pm \mathcal{J}_0$) & complex ($\mp \mathcal{J}_0$) \\
complex ($\pm \mathcal{J}_0$) & complex ($\mp \mathcal{J}_0$) & complex ($\mp \mathcal{J}_0$) & complex ($\mp \mathcal{J}_0$) \\
complex ($\pm \mathcal{J}_0$) & complex ($\mp \mathcal{J}_0$) & complex ($\pm \mathcal{J}_0$) & noncomplex \\
complex ($\pm \mathcal{J}_0$) & complex ($\mp \mathcal{J}_0$) & noncomplex  & complex ($\mp \mathcal{J}_0$) \\
complex ($\pm \mathcal{J}_0$) & noncomplex  & complex ($\pm \mathcal{J}_0$) & complex ($\pm \mathcal{J}_0$) \\
noncomplex & complex ($\pm \mathcal{J}_0$) & complex ($\pm \mathcal{J}_0$) & complex ($\pm \mathcal{J}_0$) \\
noncomplex & noncomplex & noncomplex & noncomplex
\end{tabular} 
\end{center}

\subsubsection*{Type $E_6$}
The algebra $E_6$ is represented by the Dynkin diagram
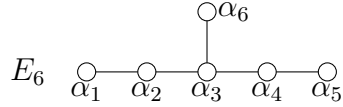
\begin{center}
  \begin{tikzpicture}[scale=.4]
    \draw (-1,0) node[anchor=east]  {$E_6$};
    \draw (0,0) circle (.3cm) node [below] {$\alpha_1$};
    \draw (0.3 cm,0) -- +(1.4 cm,0); 
    \draw (2,0) circle (.3cm) node [below] {$\alpha_2$};
    \draw (2.3 cm,0) -- +(1.4 cm,0); 
    \draw (4,0) circle (.3cm) node [below] {$\alpha_3$};    
    \draw (4.3 cm,0) -- +(1.4 cm,0); 
    \draw (6,0) circle (.3cm) node [below] {$\alpha_{4}$};
    \draw (6.3 cm,0) -- +(1.4 cm,0); 
    \draw (8,0) circle (.3cm) node [below] {$\alpha_{5}$};
    \draw (4cm,.3cm) -- +(0,1.4 cm);
    \draw (4,2) circle (.3cm) node [anchor=west] {$\alpha_6$};
  \end{tikzpicture}
\end{center}
and has $36$ positive roots. This algebra appears three times in Table \ref{tabelaflag}, but we need to analyse just two cases:
\begin{enumerate}
\item $\mathbb{F}_\Theta = E_6/U(1)\times U(1)\times SO(8)$ with $\Theta = \Sigma \backslash \{\alpha_1,\alpha_5\}$;

\item $\mathbb{F}_\Theta = E_6/SU(5)\times U(1)\times U(1)$ with $\Theta = \Sigma \backslash \{\alpha_1, \alpha_2\}$.
\end{enumerate}
because the other three cases are covered by Proposition \ref{oneroot}.

The first case, we have $\mathbb{F}_\Theta$ with $\Theta = \Sigma \backslash \{\alpha_1,\alpha_5\}$ and $\mathfrak{m} = \mathfrak{m}_1 \oplus \mathfrak{m}_2 \oplus \mathfrak{m}_3$ in which 
\begin{align*}
\mathfrak{m}_1 = & \ \{\alpha_1, \alpha_1+\alpha_2, \alpha_1+\alpha_2+\alpha_3, \alpha_1+\alpha_2+\alpha_3+\alpha_4, \alpha_1+\alpha_2+\alpha_3+\alpha_6, \\ & \alpha_1+\alpha_2+\alpha_3+\alpha_4+\alpha_6, \alpha_1+\alpha_2+2\alpha_3+\alpha_4+\alpha_6, \alpha_1+2\alpha_2+2\alpha_3+\alpha_4+\alpha_6\}
\end{align*}
\begin{align*}
\mathfrak{m}_2 = & \ \{ \alpha_5, \alpha_4+\alpha_5, \alpha_3+\alpha_4+\alpha_5, \alpha_3+\alpha_4+\alpha_5+\alpha_6, \alpha_2+\alpha_3+\alpha_4+\alpha_5, \\ & \alpha_2+\alpha_3+\alpha_4+\alpha_5+\alpha_6, \alpha_2+2\alpha_3+\alpha_4+\alpha_5+\alpha_6, \alpha_2+2\alpha_3+2\alpha_4+\alpha_5+\alpha_6\}
\end{align*}
\begin{align*}
\mathfrak{m}_3 = & \ \{ \alpha_1+\alpha_2+\alpha_3+\alpha_4+\alpha_5, \alpha_1+\alpha_2+\alpha_3+\alpha_4+\alpha_5+\alpha_6, \\ & \alpha_1+\alpha_2+2\alpha_3+\alpha_4+\alpha_5+\alpha_6, \alpha_1+\alpha_2+2\alpha_3+2\alpha_4+\alpha_5+\alpha_6, \\ & \alpha_1+2\alpha_2+2\alpha_3+\alpha_4+\alpha_5+\alpha_6, \alpha_1+2\alpha_2+2\alpha_3+2\alpha_4+\alpha_5+\alpha_6, \\ & \alpha_1+2\alpha_2+3\alpha_3+2\alpha_4+\alpha_5+\alpha_6, \alpha_1+2\alpha_2+3\alpha_3+2\alpha_4+\alpha_5+2\alpha_6\}. 
\end{align*}
It is easy to see that given $\alpha$ and $\beta$ roots such that $\alpha+ \beta$ is also a root, then $\alpha \in \mathfrak{m}_1$, $\beta \in \mathfrak{m}_2$ and $\alpha+\beta \in \mathfrak{m}_3$. Therefore, an invariant generalized almost complex structure $\mathcal{J}$ is integrable if and only if is one of the following
\begin{center}
\begin{tabular}{c|c|c}
$\mathcal{J}|_{\mathfrak{m}_1 \oplus \mathfrak{m}_1 ^\ast}$ & $\mathcal{J}|_{\mathfrak{m}_2 \oplus \mathfrak{m}_2 ^\ast}$ & $\mathcal{J}|_{\mathfrak{m}_3 \oplus \mathfrak{m}_3 ^\ast}$ \\ \hline
complex ($\pm \mathcal{J}_0$) & complex ($\pm \mathcal{J}_0$) & complex ($\pm \mathcal{J}_0$) \\
complex ($\pm \mathcal{J}_0$) & complex ($ \mp \mathcal{J}_0$) & complex ($\pm \mathcal{J}_0$) \\
complex ($\pm \mathcal{J}_0$) & complex ($\mp \mathcal{J}_0$) & complex ($\mp \mathcal{J}_0$) \\
complex ($\pm \mathcal{J}_0$) & noncomplex & complex ($\pm \mathcal{J}_0$) \\
noncomplex  & complex ($\pm \mathcal{J}_0$) & complex ($\pm \mathcal{J}_0$) \\
complex ($\pm \mathcal{J}_0$) & complex ($\mp \mathcal{J}_0$) & noncomplex \\
noncomplex & noncomplex & noncomplex
\end{tabular} 
\end{center}

The second case is given by $\mathbb{F}_\Theta$ with $\Theta = \Sigma \backslash \{ \alpha_1 ,\alpha_2\}$ where $\mathfrak{m} = \mathfrak{m}_1\oplus \mathfrak{m}_2\oplus \mathfrak{m}_3 \oplus \mathfrak{m}_4$ in which 
\[
\mathfrak{m}_1 = \{ \alpha_1\}
\]
\begin{align*}
\mathfrak{m}_2 = & \ \{ \alpha_2, \alpha_2+\alpha_3, \alpha_2+\alpha_3+\alpha_4, \alpha_2+\alpha_3+\alpha_6, \alpha_2+\alpha_3+\alpha_4+\alpha_5, \alpha_2+\alpha_3+\alpha_4+\alpha_6, \\ & \alpha_2+\alpha_3+\alpha_4+\alpha_5+\alpha_6, \alpha_2+2\alpha_3+\alpha_4+\alpha_6, \alpha_2+2\alpha_3+\alpha_4+\alpha_5+\alpha_6, \\ & \alpha_2+2\alpha_3+2\alpha_4+\alpha_5\alpha_6\}
\end{align*}
\begin{align*}
\mathfrak{m}_3 = & \ \{ \alpha_1+\alpha_2,\alpha_1+\alpha_2+\alpha_3, \alpha_1+\alpha_2+\alpha_3+\alpha_4, \alpha_1+\alpha_2+\alpha_3+\alpha_4+\alpha_5, \\ & \alpha_1+\alpha_2+\alpha_3+\alpha_6, \alpha_1+\alpha_2+\alpha_3+\alpha_4+\alpha_5+\alpha_6, \alpha_1+\alpha_2+2\alpha_3+\alpha_4+\alpha_6, \\ & \alpha_1+\alpha_2+2\alpha_3+\alpha_4+\alpha_5+\alpha_6, \alpha_1+\alpha_2+2\alpha_3+2\alpha_4+\alpha_5+\alpha_6\}
\end{align*}
\begin{align*}
\mathfrak{m}_4 = & \ \{ \alpha_1+2\alpha_2+2\alpha_3+\alpha_4+\alpha_6, \alpha_1+2\alpha_2+2\alpha_3+\alpha_4+\alpha_5+\alpha_6, \\ & \alpha_1+2\alpha_2+2\alpha_3+2\alpha_4+\alpha_5+\alpha_6, \alpha_1+2\alpha_2+3\alpha_3+2\alpha_4+\alpha_5+\alpha_6, \\ & \alpha_1+2\alpha_2+3\alpha_3+2\alpha_4+\alpha_5+2\alpha_6\}.
\end{align*}
Observe that given a triple of roots ($\alpha, \beta,\alpha+\beta$), the only options are:
\begin{itemize}
\item $\alpha \in \mathfrak{m}_1$, $\beta \in \mathfrak{m}_2$ and $\alpha+\beta \in \mathfrak{m}_3$;
\item $\alpha \in \mathfrak{m}_2$, $\beta \in \mathfrak{m}_3$ and $\alpha+\beta \in \mathfrak{m}_4$.
\end{itemize}
Therefore, an invariant integrable generalized complex structure on $\mathbb{F}_\Theta$ must be one of the following structures:
\begin{center}
\begin{tabular}{c|c|c|c}
$\mathcal{J}|_{\mathfrak{m}_1 \oplus \mathfrak{m}_1 ^\ast}$ & $\mathcal{J}|_{\mathfrak{m}_2 \oplus \mathfrak{m}_2 ^\ast}$ & $\mathcal{J}|_{\mathfrak{m}_3 \oplus \mathfrak{m}_3 ^\ast}$ & $\mathcal{J}|_{\mathfrak{m}_4 \oplus \mathfrak{m}_4 ^\ast}$  \\ \hline
complex ($\pm \mathcal{J}_0$) & complex ($\pm \mathcal{J}_0$) & complex ($\pm \mathcal{J}_0$) & complex ($\pm \mathcal{J}_0$) \\
complex ($\pm \mathcal{J}_0$) & complex ($\mp \mathcal{J}_0$) & complex ($\pm \mathcal{J}_0$) & complex ($\pm \mathcal{J}_0$) \\
complex ($\pm \mathcal{J}_0$) & complex ($\mp \mathcal{J}_0$) & complex ($\pm \mathcal{J}_0$) & complex ($\mp \mathcal{J}_0$) \\
complex ($\pm \mathcal{J}_0$) & complex ($\mp \mathcal{J}_0$) & complex ($\mp \mathcal{J}_0$) & complex ($\mp \mathcal{J}_0$) \\
complex ($\pm \mathcal{J}_0$) & complex ($\mp \mathcal{J}_0$) & complex ($\pm \mathcal{J}_0$) & noncomplex \\
complex ($\pm \mathcal{J}_0$) & complex ($\mp \mathcal{J}_0$) & noncomplex  & complex ($\mp \mathcal{J}_0$) \\
complex ($\pm \mathcal{J}_0$) & noncomplex  & complex ($\pm \mathcal{J}_0$) & complex ($\pm \mathcal{J}_0$) \\
noncomplex & complex ($\pm \mathcal{J}_0$) & complex ($\pm \mathcal{J}_0$) & complex ($\pm \mathcal{J}_0$) \\
noncomplex & noncomplex & noncomplex & noncomplex
\end{tabular} 
\end{center}

\subsubsection*{Type $E_7$}
The algebra $E_7$ is represented by the Dynkin diagram 
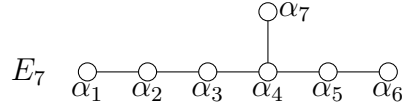
\begin{center}
  \begin{tikzpicture}[scale=.4]
    \draw (-3,0) node[anchor=east]  {$E_7$};
    \draw (-2,0) circle (.3cm) node [below] {$\alpha_1$};
    \draw (-1.7cm,0) -- +(1.4cm,0);
    \draw (0,0) circle (.3cm) node [below] {$\alpha_2$};
    \draw (0.3 cm,0) -- +(1.4 cm,0); 
    \draw (2,0) circle (.3cm) node [below] {$\alpha_3$};
    \draw (2.3 cm,0) -- +(1.4 cm,0); 
    \draw (4,0) circle (.3cm) node [below] {$\alpha_4$};    
    \draw (4.3 cm,0) -- +(1.4 cm,0); 
    \draw (6,0) circle (.3cm) node [below] {$\alpha_{5}$};
    \draw (6.3 cm,0) -- +(1.4 cm,0); 
    \draw (8,0) circle (.3cm) node [below] {$\alpha_{6}$};
    \draw (4cm,.3cm) -- +(0,1.4 cm);
    \draw (4,2) circle (.3cm) node [anchor=west] {$\alpha_7$};
  \end{tikzpicture}
\end{center}
and has $64$ positive roots. This algebra appears $4$ times in Table \ref{tabelaflag}, but we just need to analyse one case, when $\mathbb{F}_\Theta = E_7/SO(10)\times U(1)\times U(1)$ with $\Theta = \Sigma\backslash \{ \alpha_1, \alpha_2 \}$. In this case, we have $\mathfrak{m} = \mathfrak{m}_1\oplus \mathfrak{m}_2 \oplus \mathfrak{m}_3 \oplus \mathfrak{m}_4$ where 
\[
\mathfrak{m}_1 = \{ \alpha_1\}
\]
\begin{align*}
\mathfrak{m}_2 = & \ \{ \alpha_2, \alpha_2+\alpha_3, \alpha_2+2\alpha_3+2\alpha_4+2\alpha_5+\alpha_6+\alpha_7, \alpha_2+\alpha_3+2\alpha_4+2\alpha_5+\alpha_6+\alpha_7,\\ & \alpha_2+\alpha_3+\alpha_4+\alpha_5+\alpha_6+\alpha_7, \alpha_2+\alpha_3+\alpha_4+\alpha_5, \alpha_2+\alpha_3+\alpha_4, \\ & \alpha_2+2\alpha_3+3\alpha_4+2\alpha_5+\alpha_6+\alpha_7, \alpha_2+\alpha_3+2\alpha_4+\alpha_5+\alpha_6+\alpha_7, \\ & \alpha_2+2\alpha_3+2\alpha_4+\alpha_5+\alpha_6+\alpha_7, \alpha_2+2\alpha_3+3\alpha_4+2\alpha_5+\alpha_6+2\alpha_7, \\ & \alpha_2+\alpha_3+\alpha_4+\alpha_5+\alpha_6, \alpha_2+\alpha_3+\alpha_4+\alpha_5+\alpha_6+\alpha_7, \alpha_2+\alpha_3+\alpha_4+\alpha_7, \\ & \alpha_2+\alpha_3+2\alpha_4+\alpha_5+\alpha_7, \alpha_2+2\alpha_3+2\alpha_4+\alpha_5+\alpha_7\}
\end{align*}
\begin{align*}
\mathfrak{m}_3 = & \ \{ \alpha_1+\alpha_2, \alpha_1+\alpha_2+\alpha_3, \alpha_1+\alpha_2+\alpha_3+\alpha_4+\alpha_5+\alpha_6+\alpha_7, \\ & \alpha_1+\alpha_2+2\alpha_3+2\alpha_4+\alpha_5+\alpha_6+\alpha_7, \alpha_1+\alpha_2+\alpha_3+2\alpha_4+2\alpha_5+\alpha_6+\alpha_7, \\ & \alpha_1+\alpha_2+\alpha_3+\alpha_4+\alpha_5, \alpha_1+\alpha_2+2\alpha_3+3\alpha_4+2\alpha_5+\alpha_6+\alpha_7, \\ & \alpha_1+\alpha_2+\alpha_3+\alpha_4, \alpha_1+\alpha_2+2\alpha_3+2\alpha_4+\alpha_5+\alpha_6+\alpha_7,\\ & \alpha_1+\alpha_2+\alpha_3+2\alpha_4+\alpha_5+\alpha_6+\alpha_7, \alpha_1+\alpha_2+\alpha_3+\alpha_4+\alpha_5+\alpha_6, \\ & \alpha_1+\alpha_2+\alpha_3+\alpha_4+\alpha_5+\alpha_7, \alpha_1+\alpha_2+2\alpha_3+3\alpha_4+2\alpha_5+\alpha_6+\alpha_7, \\ & \alpha_1+\alpha_2+\alpha_3+\alpha_4+\alpha_7, \alpha_1+\alpha_2+2\alpha_3+2\alpha_4+\alpha_5+\alpha_7, \\ & \alpha_1+\alpha_2+\alpha_3+2\alpha_4+\alpha_5+\alpha_7\}
\end{align*}
\begin{align*}
\mathfrak{m}_4 = & \ \{ \alpha_1+2\alpha_2+2\alpha_3+2\alpha_4+2\alpha_5+\alpha_6+\alpha_7, \alpha_1+2\alpha_2+3\alpha_3+3\alpha_4+2\alpha_5+\alpha_6+\alpha_7, \\ & \alpha_1+2\alpha_2+2\alpha_3+3\alpha_4+2\alpha_5+\alpha_6+\alpha_7, \alpha_1+2\alpha_2+2\alpha_3+2\alpha_4+\alpha_5+\alpha_6+\alpha_7, \\ & \alpha_1+2\alpha_2+3\alpha_3+4\alpha_4+3\alpha_5+2\alpha_6+2\alpha_7, \alpha_1+2\alpha_2+3\alpha_3+4\alpha_4+2\alpha_5+\alpha_6+2\alpha_7, \\ & \alpha_1+2\alpha_2+3\alpha_3+3\alpha_4+2\alpha_5+\alpha_6+2\alpha_7,
\alpha_1+2\alpha_2+2\alpha_3+3\alpha_4+2\alpha_5+2\alpha_6+2\alpha_7, \\ &\alpha_1+2\alpha_2+2\alpha_3+2\alpha_4+\alpha_5+\alpha_7, \alpha_1+2\alpha_2+3\alpha_3+4\alpha_4+3\alpha_5+\alpha_6+2\alpha_7\}.
\end{align*}
Therefore, an invariant integrable generalized complex structure on $\mathbb{F}_\Theta$ must be one of the following:
\begin{center}
\begin{tabular}{c|c|c|c}
$\mathcal{J}|_{\mathfrak{m}_1 \oplus \mathfrak{m}_1 ^\ast}$ & $\mathcal{J}|_{\mathfrak{m}_2 \oplus \mathfrak{m}_2 ^\ast}$ & $\mathcal{J}|_{\mathfrak{m}_3 \oplus \mathfrak{m}_3 ^\ast}$ & $\mathcal{J}|_{\mathfrak{m}_4 \oplus \mathfrak{m}_4 ^\ast}$  \\ \hline
complex ($\pm \mathcal{J}_0$) & complex ($\pm \mathcal{J}_0$) & complex ($\pm \mathcal{J}_0$) & complex ($\pm \mathcal{J}_0$) \\
complex ($\pm \mathcal{J}_0$) & complex ($\mp \mathcal{J}_0$) & complex ($\pm \mathcal{J}_0$) & complex ($\pm \mathcal{J}_0$) \\
complex ($\pm \mathcal{J}_0$) & complex ($\mp \mathcal{J}_0$) & complex ($\pm \mathcal{J}_0$) & complex ($\mp \mathcal{J}_0$) \\
complex ($\pm \mathcal{J}_0$) & complex ($\mp \mathcal{J}_0$) & complex ($\mp \mathcal{J}_0$) & complex ($\mp \mathcal{J}_0$) \\
complex ($\pm \mathcal{J}_0$) & complex ($\mp \mathcal{J}_0$) & complex ($\pm \mathcal{J}_0$) & noncomplex \\
complex ($\pm \mathcal{J}_0$) & complex ($\mp \mathcal{J}_0$) & noncomplex  & complex ($\mp \mathcal{J}_0$) \\
complex ($\pm \mathcal{J}_0$) & noncomplex  & complex ($\pm \mathcal{J}_0$) & complex ($\pm \mathcal{J}_0$) \\
noncomplex & complex ($\pm \mathcal{J}_0$) & complex ($\pm \mathcal{J}_0$) & complex ($\pm \mathcal{J}_0$) \\
noncomplex & noncomplex & noncomplex & noncomplex
\end{tabular} 
\end{center}

Considering all the cases studied, we have covered all flag manifolds with two, three and four isotropy summands. For each case we describe all the possibilities for $\mathcal{J}$ an invariant generalized almost complex structure to be integrable. In particular, Proposition \ref{oneroot}, cover some cases with five or six isotropy summands. The same reasoning can be used to describe the invariant generalized complex structures for other flag manifolds with more than four isotropy summands.


\newpage
\begin{thebibliography}{WWW}
\bibitem[AC1]{AC1} A. Arvanitoyeorgos and I. Chrysikos: {\it Invariant Einstein Käler metrics on generalized flag manifolds with two isotropy summands}, J. Aust. Math. Soc. {\bf 90.2} (2011), 237--251.

\bibitem[AC2]{AC2} A. Arvanitoyeorgos and I. Chrysikos: {\it Invariant Einstein Käler metrics on flag manifolds with four isotropy summands}, Ann. Glob. Anal. Geom. {\bf 37.2} (2010), 185--219.

\bibitem[ACS]{ACS} A. Arvanitoyeorgos, I. Chrysikos and Y. Sakane: {\it Homogeneous Einstein metrics on generalized flag manifolds with five isotropy summands}, Int. J. Math. {\bf 24.10} (2014), 1350077.

\bibitem[C]{C} Cavalcanti, G. R., \textit{Introduction to generalized complex geometry}. IMPA, (2007).

\bibitem[K]{Kim} K. Kimura: {\it Homogeneous Einstein metrics on certain Kähler C-spaces}, Adv. Stud. Pure Math. {\bf 18-I} (1990) 303-320.

\bibitem[G1]{G1} Gualtieri, M., \textit{Generalized complex geometry}, D.Phil. thesis, Oxford University, (2003).

\bibitem[G2]{G2}   Gualtieri, M., {\it Generalized complex geometry}, Ann. of Math. (2),  {\bf 174.1} (2011), 75--123. 

\bibitem[H]{H} Hitchin, N., {\it Generalized Calabi--Yau manifolds}, Q. J. Math. \textbf{54} (2003), 281--308.

\bibitem[VS]{VS}   Varea, C. A. B.; San Martin, L. A. B.,  {\it Invariant Generalized Complex Structures on Flag Manifolds}, preprint (2018) arXiv:1810.09532v2. 
\end{thebibliography}
\end{document}